\documentclass[12pt,a4paper]{article}

\usepackage{makeidx}         
\usepackage{graphicx}        
\usepackage{psfrag}                             
\usepackage{multicol}        
\usepackage{amsmath}
\usepackage{amssymb}
\usepackage{amsfonts}
\usepackage{latexsym}
\usepackage[T1]{fontenc}
\usepackage{subfigure}                      
\usepackage{float}
\usepackage{mathrsfs}
\usepackage{mathptmx}
\usepackage{upref}      
\usepackage{hyperref}
\usepackage{color,amsthm,amsgen}
\usepackage{url}

\newtheorem{theorem}{Theorem}
\newtheorem{lemma}[theorem]{Lemma}

\numberwithin{equation}{section}

\let\l\lambda

\newcommand{\lie}{\mathop{\displaystyle
{\stackrel{\scriptstyle{\varepsilon \to 0}}{\hbox to 40 pt
{\rightarrowfill}}  } }}
\newcommand{\lik}{\mathop{\displaystyle
{\stackrel{\scriptstyle{{\rm k} \to \infty}}{\hbox to 40 pt
{\rightarrowfill}}  } }}

\begin{document}

\title{\bf Longitudinal Oscillations for Eigenfunctions in Rod Like Structures}

\author{Pablo Benavent-Ocejo\footnote{\url{pablo.benavent@alumnos.unican.es}}, \, Delfina G\'omez\footnote{\url{gomezdel@unican.es}} \, Mar\'ia-Eugenia P\'erez-Mart\'inez\footnote{\url{meperez@unican.es}} \\ \\
Universidad de Cantabria, Santander, Spain}
\date{}

\maketitle

\renewcommand{\b}{\beta}
\renewcommand{\c}{\gamma}
\newcommand{\n}{\nu}
\newcommand{\x}{\xi}
\renewcommand{\l}{\lambda}
\newcommand{\m}{\mu}
\newcommand{\h}{\chi}
\newcommand{\s}{\sigma}
\renewcommand{\t}{\theta}
\newcommand{\trans}{^{\rm T}}
\newcommand{\ds}{{\partial S}}
\renewcommand{\o}{\omega}
\newcommand{\bb}{\mathcal B}
\newcommand{\rr}{\mathcal R}
\renewcommand{\ll}{\mathcal L}
\renewcommand{\ss}{\mathcal S}
\renewcommand{\gg}{\mathcal G}
\newcommand{\vv}{\mathcal V}
\newcommand{\ww}{\mathcal W}
\newcommand{\Dminus}{(D$^{\omega-}$)}
\newcommand{\Nminus}{(N$^{\omega-}$)}
\newcommand{\Rminus}{(R$^{\omega-}$)}
\newcommand{\Vminus}{{\mathcal V}^{\omega-}}
\renewcommand{\Re}{\mathop\textrm{Re}}
\renewcommand{\Im}{\mathop{\rm Im}}
\newcommand{\mxy}{|x-y|}
\def\Wminus{{\mathcal W}^{\omega-}}
\def\regular#1{C^2(S^{#1})\cap C^1(\bar{S}^{#1})}
\def\holder#1{C^{\,\,#1,\alpha}(\partial S)}
\def\ppd#1{\dfrac{\partial}{\partial#1}}
\newcommand{\e}{\varepsilon}
\newcommand{\ee}{\varepsilon}
\newcommand{\supsig}{^{(\sigma)}}
\newcommand{\up}{\upsilon}
\newcommand{\Up}{\Upsilon}
\newcommand{\p}{\psi}
\renewcommand{\a}{\alpha}
\newcommand{\f}{\varphi}
\newcommand{\F}{\Phi}
\renewcommand{\P}{\Psi}
\renewcommand{\d}{\delta}
\renewcommand{\k}{\kappa}
\newcommand{\sums}{\sum_{\s=1}^2}
\newcommand{\summ}{\sum_{m=0}^\infty}
\newcommand{\intl}{\int\limits}
\newcommand{\pd}{\partial}

\newcommand{\Cc}{{\mathcal C}}
\newcommand{\Dd}{{\mathcal D}}
\newcommand{\Fill}{\makebox[6pt]{\hfill}}

\newcommand{\tr}{{\mathrm T}}

\allowdisplaybreaks

\abstract{We  consider the spectrum of the Laplace operator on  3D rod structures, with a small cross section depending on a small parameter $\ee$. The boundary conditions are of Dirichlet type on the basis of this structure and Neumann on the lateral boundary. We focus on the low frequencies. We study the asymptotic behavior of  the eigenvalues and associated eigenfunctions, which are approached as $\ee\to 0$ by those of a 1D model with Dirichlet boundary conditions, but which takes into account the geometry of the domain. Explicit and numerical computations enlighten the interest of this study, when the parameter becomes smaller. At the same time they show that in order to capture oscillations in the transverse direction we need to deal with the high frequencies. For prism like domains, we show the different asymptotic behavior of the spectrum depending on the boundary conditions.}

\section{Introduction}\label{sec:ben1}

In this paper, we address the asymptotic behavior of the eigenvalues and eigenfunctions  for the Laplacian in  thin rod structures,  when the diameter of the transverse section tends to zero. Namely,    3D  domains of size $O(1)$ along the longitudinal direction and  $O(\ee)$ in the two other directions which are referred to as {\em transverse directions}.

These kinds of structures appear in many engineering constructions or engineering devices containing  bar or thin tube structures, but  from the  mathematical viewpoint there are clear gaps in the   description of eigenvalues and eigenfunctions  for 3D structures, in their dependence of the small parameter; specially the model with the mixed boundary conditions that we consider here is still an open problem in the literature of Applied Mathematics. The interest from the dynamical viewpoint is evident, both for models arising in diffusion  or vibrations of tube structures and / or  multistructures (cf. \cite{PanasenkoPerez2007,AmosovGomezPanasenkoPerez2024,Panasenko_book}).

Explicit computations on particular geometries of these structures (prims-like) show the different asymptotic  behavior of the eigenvalues and eigenfunctions depending on the boundary conditions. They are important to  enlighten  the order of magnitude of  the low frequencies, its asymptotic behavior, as $\ee\to 0$, and the behavior of the associated eigenfunctions. These explicit computations also show that in order to capture oscillations of the eigenfunctions different from transversal ones, for Dirichlet boundary conditions,  or longitudinal ones in the rest of cases we need to deal with the high frequencies. Also they show how different the above-mentioned behavior is depending on the boundary conditions.

Although for Dirichlet (Neumann respectively)  boundary conditions the behavior of the spectra has been outlined in \cite{CardoneDranteNazarov2010}   (\cite{BorisovCardone2011,AmosovGomezPanasenkoPerez2024,Panasenko_book}, respectively), the structure of the associated eigenfunctions needs to be clarified. This is why in Section \ref{sec:ben3} we provide explicit computations along with graphics illustrating different phenomena for the eigenfunctions. The figures have been obtained by means of numerical computations using the PDE Toolbox of Matlab 2024b and show how important an asymptotic analysis is in the case where the explicit computations do not work. As a matter of fact, we notice numerical instabilities  when the diameter $O(\varepsilon)$ becomes smaller (already $\varepsilon =0.01$ provides such instabilities).

Also, experimentally, we have detected important differences between 2 and 3 dimensions, since, for instance, a small perturbation along the longitudinal direction for Dirichlet boundary conditions does not imply a localization of the eigenfunctions associated to the low frequencies, contrarily to what one can detect for the 2D rod structure (cf. \cite{FriedlanderSolomyak} and Figure \ref{fig:ben1}).

We need a thorough study of the asymptotic behavior of the eigenelements.  We focus on the low frequencies for  mixed boundary conditions  and  provide a 1D limit model with Dirichlet conditions, which takes into account the geometry of the domain and  approaches  the spectrum of the original problem (cf. \eqref{eq:ben24}).

For thin 2D rod structures and 3D like films structures with only one of the dimensions smaller than the other or an oscillating boundary,   we refer to \cite{NazarovPerezTaskinen2016, ArrietaNakasatoVillanueva2025}  and references therein. The junction of rod structures or  thin films has been addressed in \cite{GaudielloSili2007, GaudielloGomezPerez2023, Chesnel}.

The structure of the paper is as follows. Section \ref{sec:ben2} contains the statement of the problem under consideration. Section \ref{sec:ben3} provides explicit computations for a prism and illustrates the phenomena by means of numerical computations. Section \ref{sec:ben4} contains the limit problem and the statement of the main result of convergence. Finally, in the appendix, cf. Section \ref{sec:ben5}, we present illustrative computations for the case of Dirichlet Laplacian or Neumann Laplacian.

\section{The statement of the problem}\label{sec:ben2}

Let $G$ be an open bounded domain of $ \mathbb{R}^3$ with a Lipschitz boundary,  which for the sake of simplicity we assume to be placed along the $x_1$ axis (cf. \eqref{eq:ben1}). Let $\partial G$ denote the boundary of $G$ which is assumed to be the union  two plane faces, which we denote by $\Gamma_0$ and $\Gamma_1$, and another  lateral surface $\Gamma_l$ of   $ \mathbb{R}^3$. In particular,   we set $\overline{\Gamma_0}=  \overline G\cap \{x_1=\ell_0\}$, $\overline{ \Gamma_1}= \overline G\cap \{  x_1=\ell_1 \}$, for any fixed  positive constants $\ell_0\leq 0 < \ell_1$ and
$$\partial G=\overline{\Gamma_0}\cup \overline{\Gamma_1}\cup  \overline{\Gamma_l}.$$
Namely, $G$ is a like prism domain of the space, with  a somehow arbitrary lateral surface ${\Gamma_l}$, that admits a representation
\begin{equation}\label{eq:ben1}
G=\bigcup_ {x_1 \in (\ell_0,\ell_1)} \{(x_1, x_2,x_3)\, :\, (x_2,x_3)\in D_{x_1}\},
\end{equation}
$D_{x_1}$ being the transverse sections;  namely, for any fixed $x_1\in  (\ell_0,\ell_1)$, $D_{x_1}$  is  an open  domain of the plane with a Lispchitz boundary which depends on $x_1$ (cf., for example, \eqref{eq:ben9}, \eqref{eq:ben10}, \eqref{eq:ben12} and \eqref{eq:ben14}). In this way, when  $D_{x_1}=D,\, \forall x_1 \in (\ell_0,\ell_1) $, $G$ is a tube domain
$G= (\ell_0,\ell_1)\times D$.
In addition, we assume that the area of the cross sections $D_{x_1}$ satisfies,
\begin{equation}\label{eq:ben2}
0<c_0<\vert D_{x_1}\vert \leq c_1, \quad \forall x_1\in [\ell_0,\ell_1],
\end{equation}
for certain constants $c_0$ and $c_1$ independent of $x_1$.

Let $\ee$ denote a small parameter $\ee\in (0,1)$ that we shall make to go to 0. We consider $G_\ee$ to be the domain
\begin{equation*}
G_\ee:=\{(x_1,x_2,x_3) \, :\, (x_1,\, \ee^{-1} x_2 ,\,\ee^{-1} x_3 )\in G\} \, .
\end{equation*}
Let $\Gamma_\ee^D$ denote the two faces perpendicular to the $x_1$ axis,   namely:
$$\Gamma_\ee^D=  \Gamma_0^\ee   \cup  \Gamma_1^\ee  \quad \mbox{ with }   \quad  \overline{ \Gamma_0^\ee} =\overline{G_\ee}\cap \{x_1=\ell_0\}\quad \mbox{and} \quad \overline{\Gamma_1^\ee } =\overline{G_\ee}\cap \{x_1=\ell_1\}.$$
Finally, the lateral surface reads $$\Gamma^l_\ee=\partial G_\ee \setminus \Gamma_\ee^D,$$
and $G_\ee$ has the representation
$$G_\ee=\bigcup_ {x_1 \in (\ell_0,\ell_1)} \{(x_1, x_2,x_3)\, :\, (x_2,x_3)\in D^\ee_{x_1}\},  \quad \mbox{ with }
D^\ee_{x_1}=\ee D_{x_1}. $$

 \smallskip

In $G_\ee$ we  consider the following eigenvalue problem with mixed boundary conditions
\begin{equation}\label{eq:ben4}
\left\{\begin{array}{ll}-\Delta u^\ee=\lambda^\ee u^\ee \hbox{ in }G_\ee, \vspace{0.2cm}\\
u^\ee=0 \hbox{ on } \Gamma_\ee^D , \vspace{0.2cm}\\
\displaystyle{ \frac{\partial u^\ee}{\partial\nu}=0 \hbox{ on } \partial
G_\ee\setminus\Gamma_\ee^D,}
\end{array}\right.
\end{equation}
where   $\nu$  denotes the outward
 unit normal to $G_\ee$, $\lambda^\ee$ stands for the eigenvalue with corresponding eigenfunction $u^\ee$.

 The weak formulation of \eqref{eq:ben4}  reads: find $(\lambda^\ee, u^\ee)\in \mathbb{R}  \times H^1 (G_\ee,\Gamma_\ee^D)$, $u^\ee\not\equiv 0$,  satisfying
\begin{equation}\label{eq:ben5}
 \int_{G_\ee} \nabla u^\ee.\nabla v  \, dx=\lambda^\ee  \int_{G_\ee} u^\ee v \,dx,  \quad\forall\, v\in H^1 (G_\ee,\Gamma_\ee^D),
\end{equation}
where $H^1 (G_\ee,\Gamma_\ee^D)$ denotes the  space completion of
$$\{ u\in {\cal C}^\infty(\overline G_\ee)\,:\, u=0 \mbox{ on } \Gamma_\ee^D\},$$
equipped      with the norm generated by the scalar product $$(\nabla u,\nabla v)_{L^2(G_\ee)}.$$
Note that on account of the Poincar\'e inequality, this norm is equivalent to the usual one  in $H^1(G_\ee)$.

The formulation \eqref{eq:ben5} is classical in the couple of Hilbert spaces $H^1 (\Omega,\Gamma_\ee^D)\subset L^2(G_\ee)$ with a dense and compact embedding (cf. e.g. Section I.5 of \cite{SaSa89}), and therefore  the problem has a discrete spectrum.

For each fixed $\ee>0$,  let us denote   by
$$ 0<\lambda_1^\ee \le\lambda_2^\ee \le\cdots \lambda_n^\ee \le \dots\to \infty, \quad \mbox{ as } n\to \infty , $$ the increasing sequence of eigenvalues, where we have adopted the convention repeated eigenvalues according to their multiplicities.  Also, we consider the corresponding set  of eigenfunctions $\{u_n^\ee\}_{ n=1}^\infty$ that  can be chosen to form an orthogonal base in  $H^1 (G_\ee,\Gamma_\ee^D)$ and     in $L^2(G_\ee)$, subject to the normalization condition
\begin{equation}\label{eq:ben6}
\int_{G_\ee} \vert u^\ee\vert^2\, dx =\ee^2,\quad
\mbox{ or, equivalently, }  \quad
\int_{G } \vert u^\ee\vert^2\, dy =1,
\end{equation}
where $y $ denotes an auxiliary variable, the so-called {\em   stretching variable}.
Its connection with $x$ is given by a change of variable  which transforms $G_\ee$ into $G$, namely,
\begin{equation}\label{eq:ben7}
y_1= x_1, \quad y_2=\frac{x_2}{\ee}, \quad y_3=\frac{x_3}{\ee}.
\end{equation}

\smallskip

First,  based on the minimax principle, we show that for fixed $\ee$ the eigenvalues $\lambda_n^\ee$ are bounded from below and from above as stated in the following lemma:

\begin{lemma}\label{lemma:ben1}
Let us assume  the hypothesis of uniform boundedness \eqref{eq:ben2}.
Then, for each $n\in  \mathbb{N}$, we have the uniform bound:
\begin{equation}\label{eq:ben8} 0 < C\leq \lambda_n^\ee\leq C_n  \, \quad \quad  \forall \ee>0,\end{equation}
where    $ C $ and $C_n$  are constants independent of $\ee$.
\end{lemma}
\begin{proof}
First, for smooth functions $v$ vanishing on $\Gamma_\ee^D$, and $(x_1,x_2,x_3)\in G_\ee$, we perform an integration by parts and apply Cauchy-Schwartz inequality to get the inequality
$$\vert v(x_1,x_2,x_3)\vert^2 \leq C \int_{\ell_0}^{x_1} \vert\partial_{x_1}  v (t,x_2,x_3)\vert^2  dt $$
for $C=(\ell_0-\ell_1)^2$, and consequently, taking integrals  for  $(x_2,x_3)\in D_{x_1}^\ee$ and then, for $x_1\in (\ell_0,\ell_1)$, we have  $$\int_{G_\ee} v^2\, dx   \leq C  \int_{G_\ee} \vert\nabla  v \vert^2\,dx, $$
 where $C$ is a constant independent of $\ee$ and $v$. Using    a density argument, we deduce the Poincar\'e inequality:
 $$\int_{G_\ee} v^2\, dx   \leq C  \int_{G_\ee} \vert\nabla  v \vert^2\,dx,
 \quad \forall v\in H^1 (G_\ee,\Gamma_\ee^D),  $$
 with the constant $C$ independent of $\ee$.
Therefore,  the  left-hand side inequality of \eqref{eq:ben8} holds true.

    Now, let us show the right-hand side inequality.
  Because of the minimax principle,  we write
  \begin{equation}\nonumber
 \lambda_n^\ee =
 \min_{ E_n  \subset H^1 (G_\ee,\Gamma_\ee^D) \,\,} \max_{v \in E_n , v  \not\equiv 0}
\frac{\int_{ G_\ee} \vert\nabla v\vert^2   \, dx   }{\int_{G_\ee  } v ^2 dx},\end{equation}
where the minimum has been taken over the set of all the subspaces $E_n$ of $H^1 (G_\ee,\Gamma_\ee^D)$ of dimension $n$.

Prescribing that $u(x_1,x_2,x_3):=u(x_1)$ for each $u\in C_0^\infty(\ell_0,\ell_1)$, we can take
    the particular space $E_n^*$ generated by the   eigenvectors $[u_1^0 , u_2^0, \cdots, u_n^0]$
  corresponding to the eigenvalues $\{\lambda_1^0, \lambda_2^0, \cdots \lambda_n^0\}$  of the Dirichlet problem in $(\ell_0,\ell_1)$:
$$u''+\lambda^0 u=0, \quad x_1 \in (\ell_0,\ell_1),\quad  u(\ell_0)=u(\ell_1)=0. $$
Then, we have
  \begin{equation}\nonumber
 \lambda_n^\ee \leq  \max_{v \in E_n^* , v  \not\equiv 0}
\displaystyle \frac{\displaystyle \int_{\ell_0}^{\ell_1} \int_{D_{x_1}}  (v'(x_1))^2 \, dx_1 dx_2 dx_3     }{\displaystyle \int_{\ell_0}^{\ell_1} \int_{D_{x_1}}     (v (x_1))^2\, dx_1  dx_2 dx_3  } \leq \frac{c_1}{c_0}  \,
  \max_{v \in E_n^* , v  \not\equiv 0}
\displaystyle \frac{\displaystyle\int_{\ell_0}^{\ell_1} (v'(x_1))^2 \, dx_1      }{\displaystyle\int_{\ell_0}^{\ell_1} (v (x_1))^2\, dx_1 } =\frac{c_1}{c_0}  \lambda_n^0   .
\end{equation}

Therefore, the lemma is proved.
\end{proof}

The aim of this paper is to detect the asymptotic behavior the eigenvalues and eigenfunctions of \eqref{eq:ben5} as $\ee \to 0$. Estimates \eqref{eq:ben8} and the normalization condition \eqref{eq:ben6} ensure that, at least by subsequences of $\ee$, still denoted by $\ee$,  we can  find converging sequences of eigenvalues and eigenfunctions $(\lambda_n^\ee, u_n^\ee)$ in a suitable space. In Section \ref{sec:ben4} we identify the limit problem and state the main result of convergence.
For the sake of completeness, below we introduce different geometrical configurations  for $G_\ee$, where the convergence results work.

\subsection{Some geometrical configurations}\label{sec:ben2_1}
For the sake of completeness we introduce here some domains $G$   satisfying all the hypothesis of the section included that of Lemma \ref{lemma:ben1} with  $\ell_1>0$ fixed.

\begin{itemize}
\item The prism
\begin{equation}\label{eq:ben9}
G =(0,\ell_1)\times (0,1)\times (0,1) \quad \mbox{ while }\quad
G_\ee=(0,\ell_1)\times (0,\ee)\times (0,\ee).
\end{equation}
\item The domain composed by the union of two prisms
\begin{equation}\label{eq:ben10}
G=(0,2^{-1}\ell_1)\times D \,\, \cup \, \,
\{2^{-1}\ell_1\}\times 2^{-1} D \,\, \cup \,\,  (2^{-1}\ell_1,\ell_1) \times 2^{-1}D
\end{equation}
with $D$ the square $D=(-1,1)\times (-1,1)$, while
\begin{equation}\label{eq:ben11}
G_\ee=(0,2^{-1}\ell_1)\times \ee D \,\, \cup \, \,
\{2^{-1}\ell_1\}\times \ee 2^{-1} D \,\, \cup \,\,  (2^{-1}\ell_1,\ell_1) \times \ee 2^{-1}D.\end{equation}
\item The domain
\begin{equation}\label{eq:ben12}
G =\{(y_1,y_2,y_3)\, : \, y_1\in  (0,\ell_1), y_2\in (0,1),\, y_3\in (0, h(y_1))\}
\end{equation}   while
$   G_\ee =\{(x_1,x_2,x_3)\, : x_1\in (0,\ell_1),\, x_2\in (0,\ee), \, x_3\in (0, \ee h(x_1))\}$,
where $h$ is a Lipschitz function uniformly bounded, namely
\begin{equation}\label{eq:ben13}
0<c_0<h(x_1)<c_1, \quad \forall x_1\in [0,\ell_1]
\end{equation}
for certain constants $c_0$ and $c_1$ independent of $x_1$.
\item Similarly, we can take
\begin{equation}\label{eq:ben14}
G =\{(y_1,y_2,y_3)\,:\, y_1\in (0,\ell_1),\, y_2\in (0, h(y_1)),\, y_3\in   (0,1)\}\end{equation}   while   $ G_\ee =\{(x_1,x_2,x_3)\,:\, x_1\in (0,\ell_1),\,  x_2\in (0, \ee h(x_1)),\, x_3\in (0,\ee)\}$.
\item Also, more general domains such as
$$
G =\{(y_1,y_2,y_3)\,:\, y_1\in  (\ell_0,\ell_1) ,\, y_2\in (-h_1(y_1), h_2(y_1)) ,\, y_3\in  (-h_3(y_1), h_4(y_1)) \}  $$
while
\begin{equation}\label{eq:ben15}
\begin{split}G_\ee =&\{(x_1,x_2,x_3)\,:\, x_1\in (\ell_0,\ell_1),\, x_2\in   (-\ee h_1(x_1),\ee h_2(x_1)), \vspace{0.25cm} \\
&\quad  x_3\in (-\ee h_3(x_1),\ee h_4(x_1))\},
\end{split}
\end{equation}
with the functions $h_j$ satisfying \eqref{eq:ben13} while $j=1,2,3,4$, and always   provided that $G$ is an open domain with a Lipschitz boundary.
\end{itemize}

\section{Explicit and numerical computations} \label{sec:ben3}

When the domain $G$ is the prism $G =(0,\ell_1)\times (0,1)\times (0,1)$ and $G_\ee=(0,\ell_1)\times (0,\ee)\times (0,\ee)$,
using separation of variables, we can compute explicitly the eigenvalues and the corresponding eigenfunctions of \eqref{eq:ben4}. Indeed, if we look for $u^\ee$ in the form $u^\ee(x_1,x_2,x_3)=F(x_1)G(x_2)H(x_3)$ for certain functions $F,\, G$ and $H$, we get to the following 1D spectral problems
\begin{eqnarray}
&F''(x_1)+\mu_1 F(x_1)=0  \quad x_1 \in (0,\ell_1),\quad  F(0)=F(\ell_1)=0, \vspace{0.25cm} \label{eq:ben16}\\
&G''(x_2)+\mu_2^\ee F(x_2)=0  \quad x_2 \in (0,\ee),\quad  G'(0)=G'(\ee)=0, \vspace{0.25cm} \label{eq:ben17}\\
&H''(x_3)+\mu_3^\ee H(x_3)=0  \quad x_3 \in (0,\ee),\quad  H'(0)=H'(\ee)=0, \label{eq:ben18}
\end{eqnarray}
where $\lambda^\ee=\mu_1+\mu_2^\ee+\mu_3^\ee$, $\mu_1, \mu_2^\ee, \mu_3^\ee\in\mathbb{R}$  and $F\not\equiv 0, \, G\not\equiv 0, \, H\not\equiv 0$. Solving the above problems, we obtain the eigenvalues of \eqref{eq:ben4}, which depend on three parameters and are given by
\begin{equation}\label{eq:ben19}
\lambda_{mrs}^\ee =\Big(\dfrac{m\pi}{\ell_1}\Big)^2+\Big(\dfrac{r\pi}{\ee}\Big)^2 + \Big(\dfrac{s\pi}{\ee}\Big)^2, \quad  m\in\mathbb{N},\,\, r, s \in \mathbb{N}\cup\{0\}
\end{equation}
The corresponding eigenfunctions are
\begin{equation}\label{eq:ben20}
u_{mrs}^\ee =A_{mrs}\sin\Big(\dfrac{m\pi x_1}{\ell_1}\Big)\cos\Big(\dfrac{r\pi x_2}{\ee}\Big)\cos\Big(\dfrac{s\pi x_3}{\ee}\Big),  \quad A_{mrs}\in\mathbb{R},  \,\,   m\in\mathbb{N},\,\, r, s \in \mathbb{N}\cup\{0\}.
\end{equation}

\begin{figure}[t]
\begin{center}
\scalebox{0.4}{\includegraphics{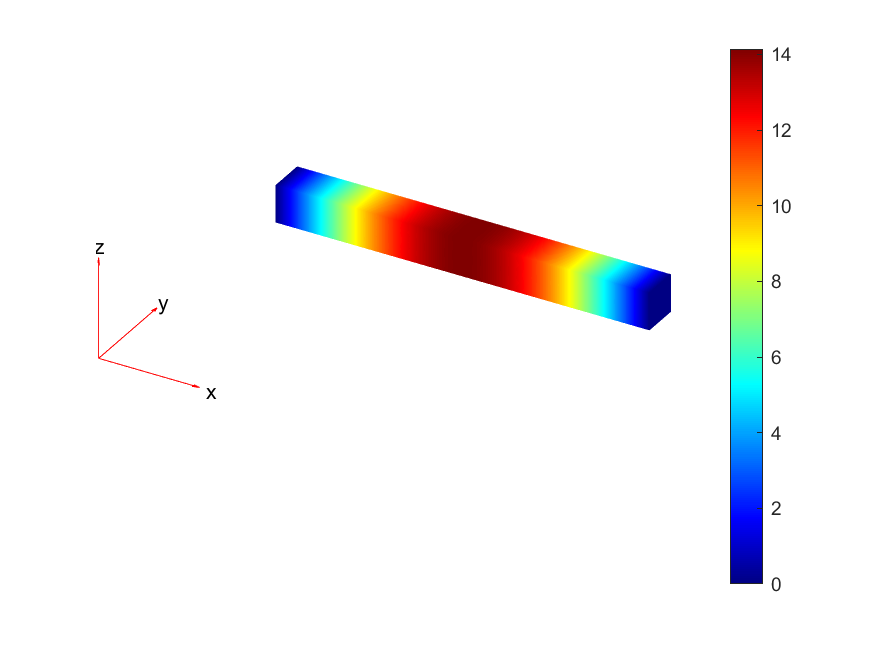}}
\scalebox{0.4}{\includegraphics{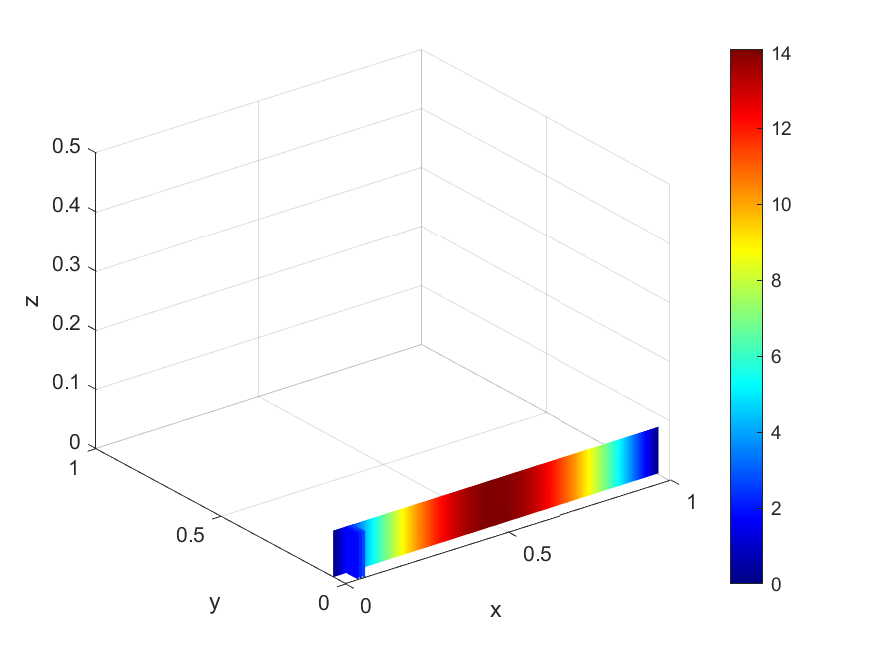}}
\scalebox{0.4}{\includegraphics{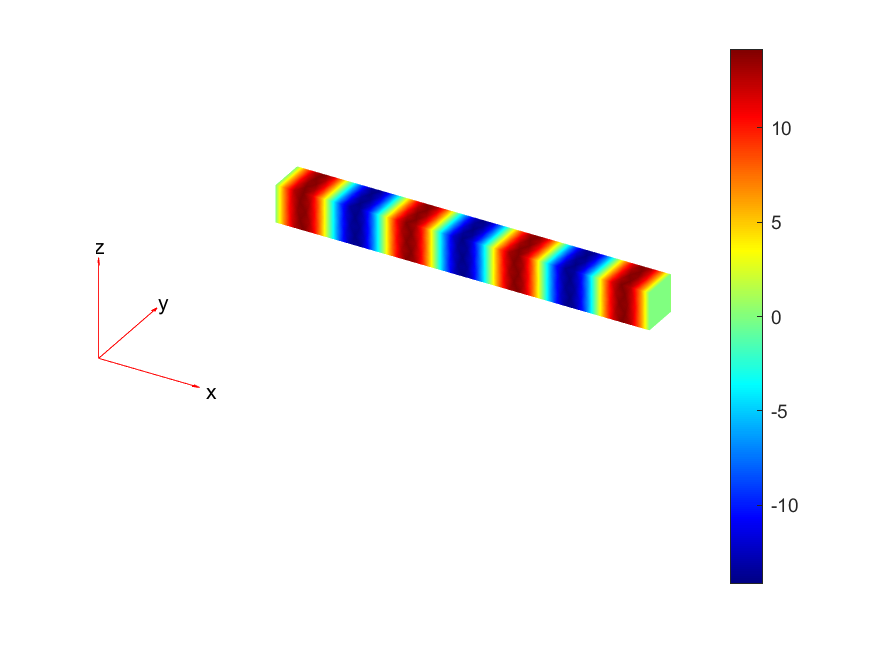}}
\scalebox{0.4}{\includegraphics{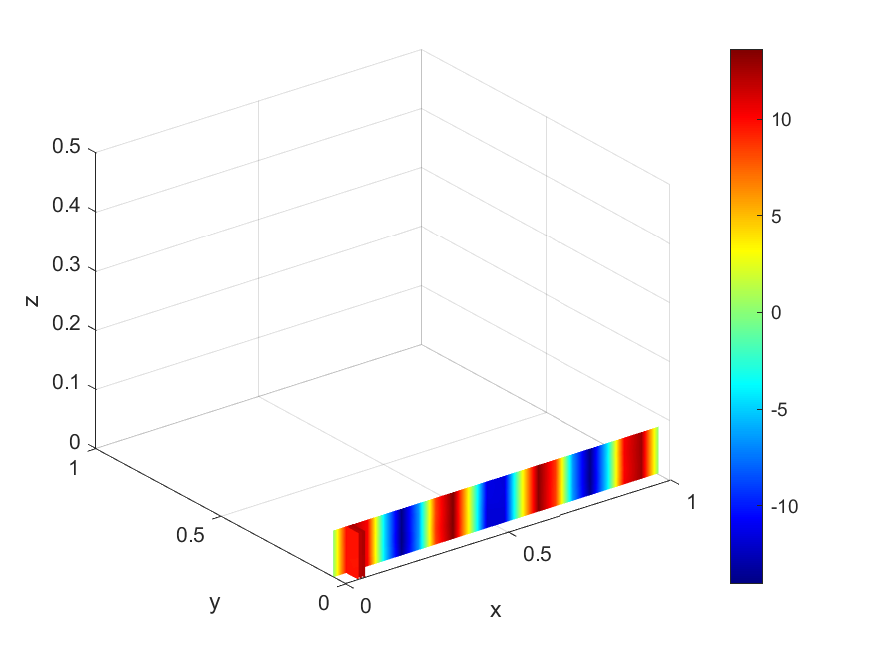}}
\scalebox{0.4}{\includegraphics{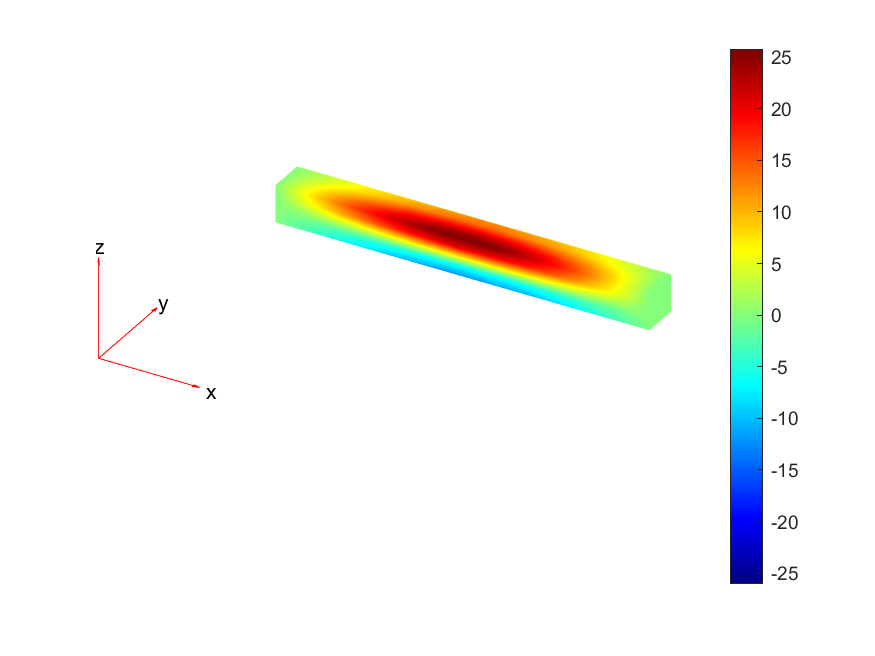}}
\scalebox{0.4}{\includegraphics{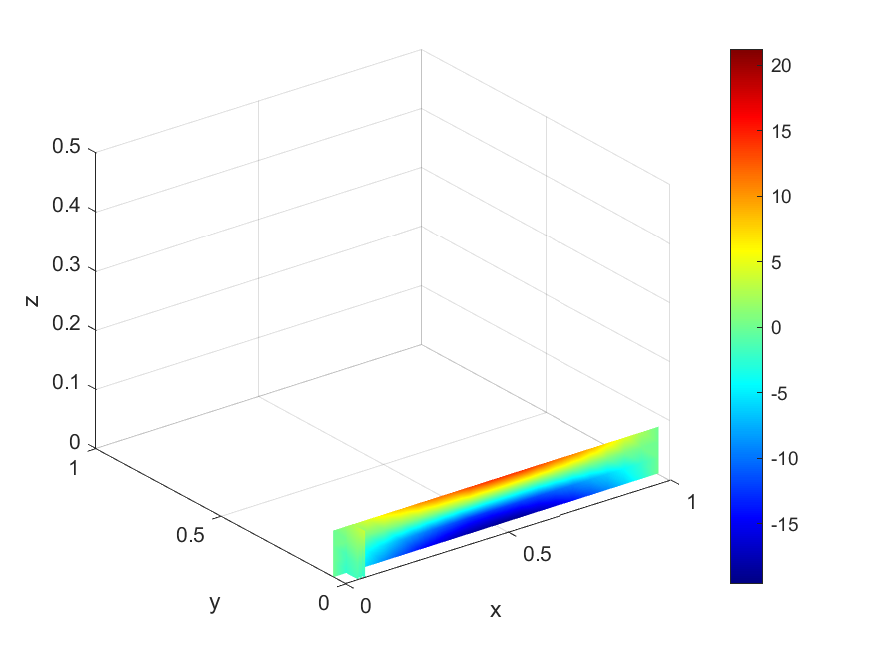}}
\end{center}
\caption{\small Approximations of eigenfunctions of \eqref{eq:ben4} with $G_\ee=(0,1)\times (0,\ee)\times (0,\ee)$ and $\ee=0.1$. The figures are obtained choosing the eigenvalues $\lambda_1^\ee=\pi^2\approx 9.87,$ $\lambda_7^\ee=49\pi^2\approx 484.08$ and $\l_{11}^\ee=101\pi^2\approx 1002.85.$ } \label{fig:ben2}
\end{figure}

From \eqref{eq:ben19} and \eqref{eq:ben20}, we observe that, for $\ee$ small, the oscillations of the eigenfunctions corresponding to the low frequencies are longitudinal (associated with the parameters $r=s=0$), that is, oscillations occur along the $x_1$ axis. In order to capture transverse  oscillations of the eigenfunctions we need to deal with the high frequencies $\lambda^\ee=O(\ee^{-2})$ (cf. Figure \ref{fig:ben2} where the last two graphics correspond to $r=0$ and  $m=s=1$).

Let us note that for other domains $G$, these computations are cumbersome or simply do not work. Therefore,  to illustrate this phenomenon we show the graphics of several eigenfunctions of problem \eqref{eq:ben4} posed in some domains described in Section \ref{sec:ben2_1},  which  have been obtained by means of numerical computations.   To do it, the domains have been generated with Blender, a free and open-source 3D computer graphics software, while for the representation of the eigenfunctions we have used the command {\em solvepdeeig} in the Partial Differential Equations Toolbox of Matlab 2024b, once we have imported the geometry of the domain and described the boundary conditions.

Figure \ref{fig:ben2} shows numerical approximations of the eigenfunctions corresponding to the first, seventh and eleventh eigenvalue of \eqref{eq:ben4} when the domain $G_\ee$ is a prism $G_\ee=(0,1)\times (0,\ee)\times (0,\ee)$ and $\ee=0.1$ (see \eqref{eq:ben9}). Here, $\lambda_1^\ee=\pi^2\approx 9.87,$ $\lambda_7^\ee=49\pi^2\approx 484.08$ and it is necessary to reach the eleventh eigenvalue $\l_{11}^\ee=101\pi^2\approx 1002.85$, which corresponds to $r=0$ and  $m=s=1$,  to capture transversal oscillations.
In the 3D figures, on the left hand side, with the color we see the oscillation of the eigenfunctions along the surface. In the plane figures, on the right hand side, we observe the oscillations  associated with the corresponding eigenfunctions along some cross sections. This graphical framework repeats for the rest of the figures.
Note that the first eigenfunction is always positive. For the rest of the eigenfunctions, in the figures,   changes from cold colors to warm colors mean oscillations in the suitable directions.

\begin{figure}[t]
\begin{center}
\scalebox{0.4}{\includegraphics{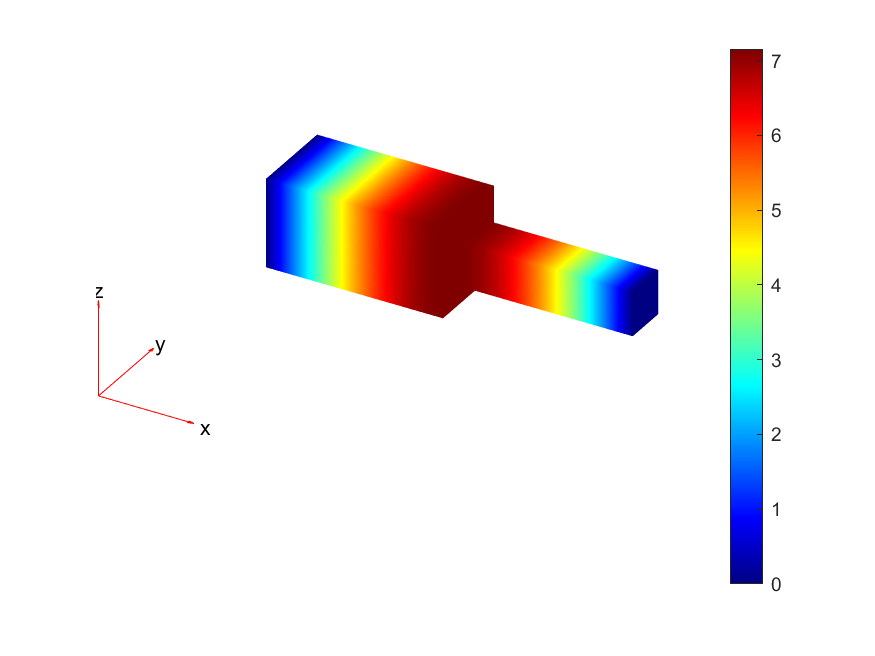}}
\scalebox{0.4}{\includegraphics{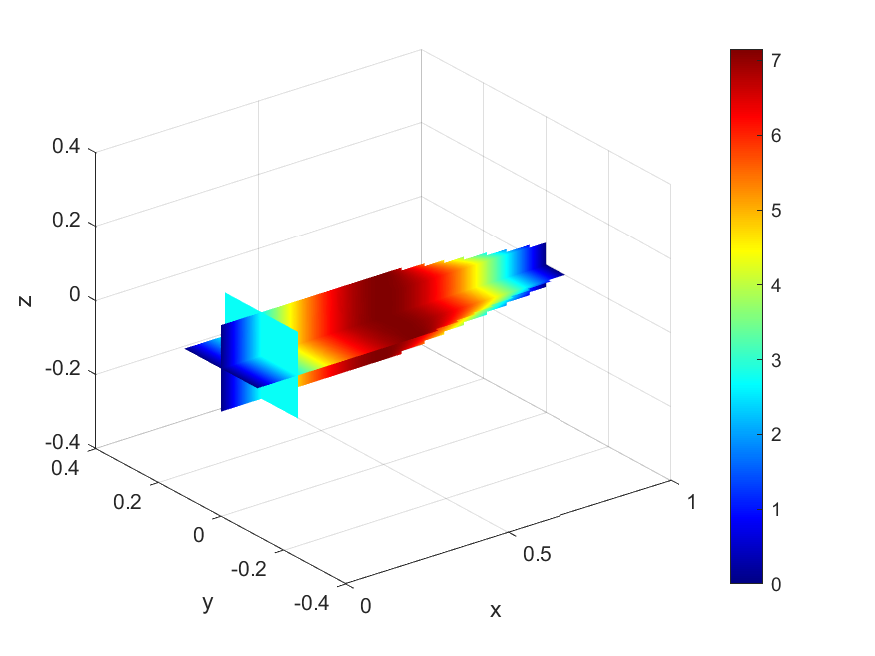}}
\scalebox{0.4}{\includegraphics{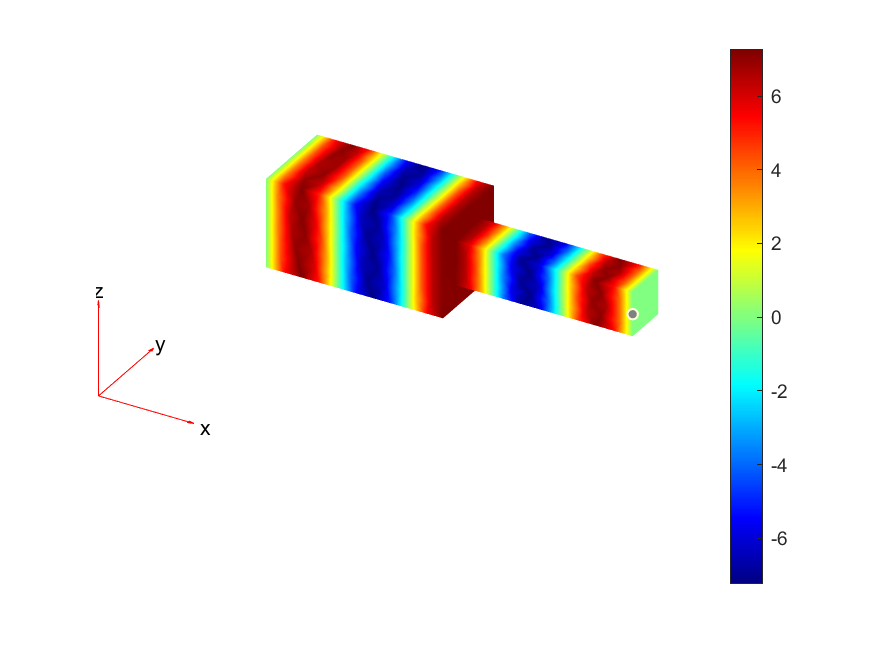}}
\scalebox{0.4}{\includegraphics{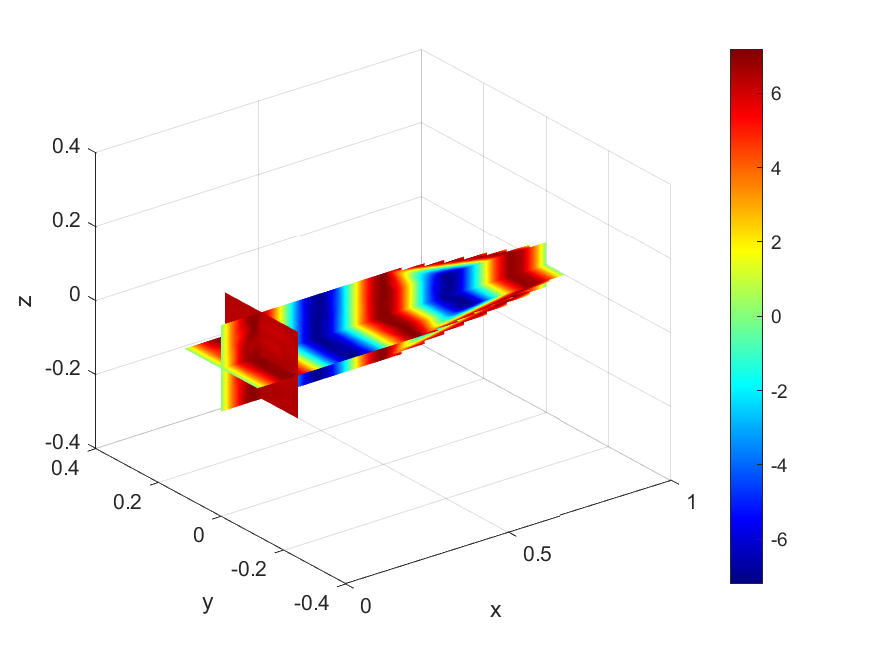}}
\scalebox{0.4}{\includegraphics{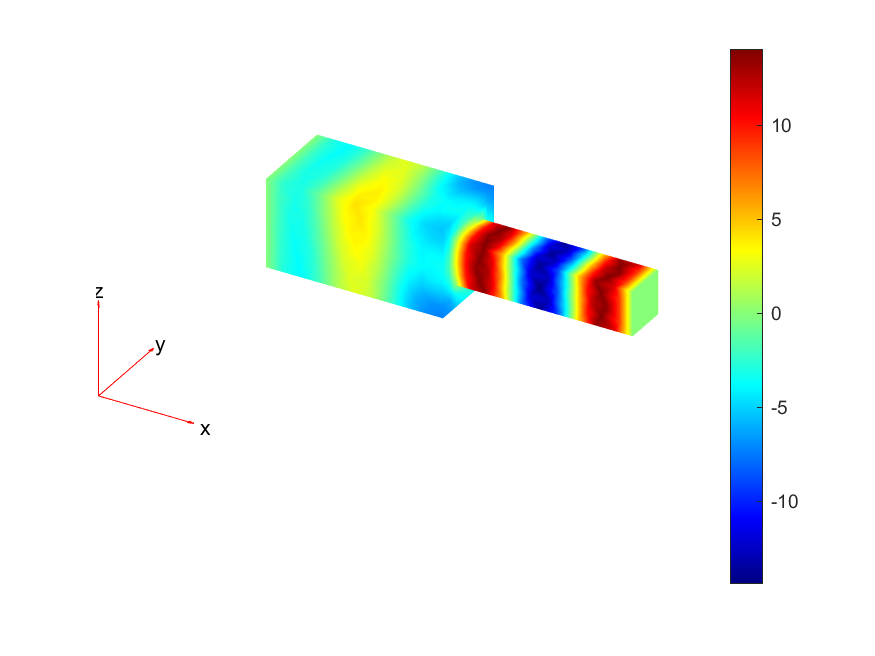}}
\scalebox{0.4}{\includegraphics{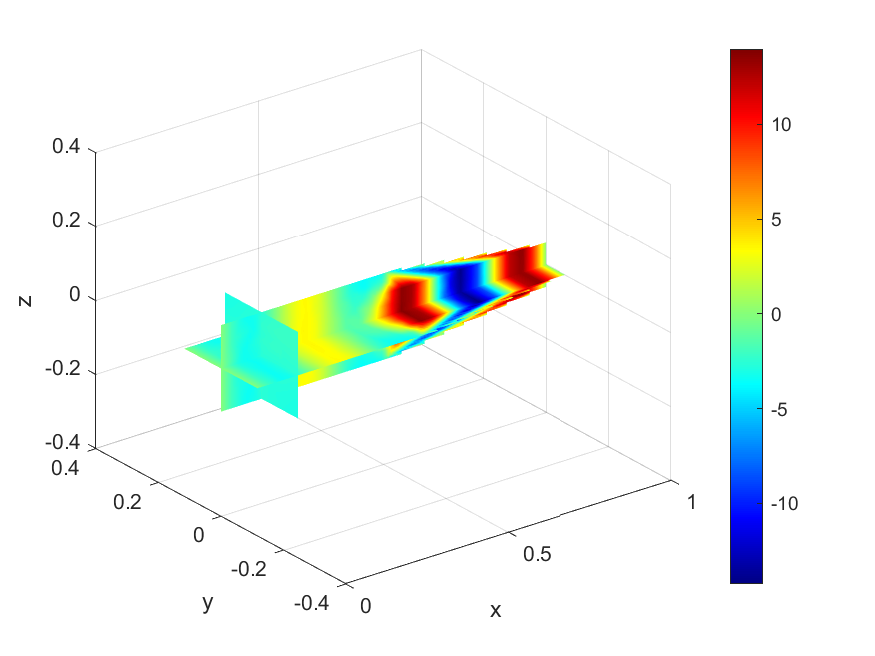}}
\end{center}
\caption{\small Approximations of eigenfunctions of \eqref{eq:ben4} with $G_\ee=(0,2^{-1})\times(-\ee,\ee)\times(-\ee,\ee) \, \cup \,
\{2^{-1}\}\times (-\ee,\ee)\times(-\ee,\ee) \, \cup \, (1/2,1)\times (-\ee2^{-1},\ee2^{-1})\times (-\ee2^{-1},\ee2^{-1})$ and $\ee=8^{-1}$. The figures are obtained choosing the eigenvalues $\lambda_1^\ee \approx 9.87,$ $\lambda_7^\ee\approx 247.50$ and $\l_{11}^\ee\approx 332.88.$} \label{fig:ben3}
\end{figure}

Figure \ref{fig:ben3} shows numerical approximations of the eigenfunctions corresponding to the first, seventh and eleventh eigenvalue of \eqref{eq:ben4} when the domain $G_\ee=(0,2^{-1})\times(-\ee,\ee)\times(-\ee,\ee) \cup
\{2^{-1}\! \}\times (-\ee 2^{-1}\! ,\ee 2^{-1}\!)\times(-\ee2^{-1}\! ,\ee2^{-1}\!) \cup (2^{-1}\!,1)\times (-\ee2^{-1}\! ,\ee2^{-1})\times (-\ee2^{-1}\!,\ee2^{-1}\!)$ and $\ee=8^{-1}$ (see \eqref{eq:ben11}).  Here, $\lambda_1^\ee \approx 9.87,$  $\lambda_7^\ee\approx 247.50$ and $\l_{11}^\ee\approx 332.88.$

\begin{figure}[t]
\begin{center}
\scalebox{0.4}{\includegraphics{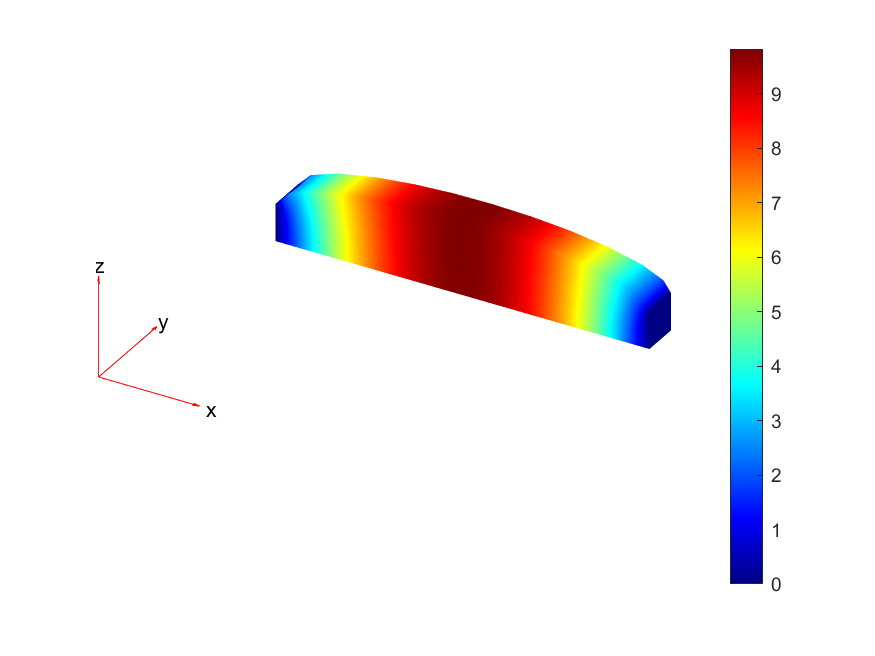}}
\scalebox{0.4}{\includegraphics{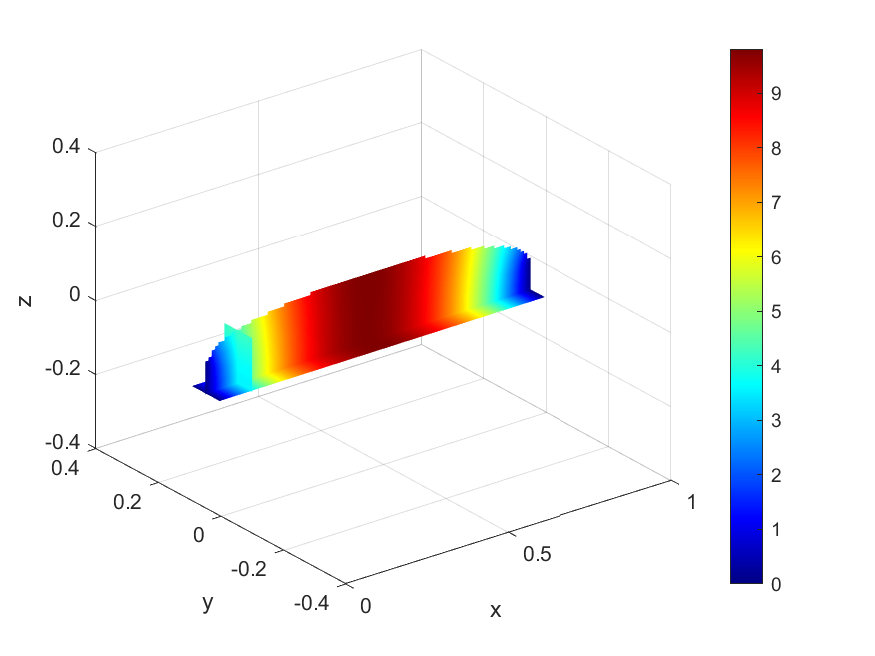}}
\scalebox{0.4}{\includegraphics{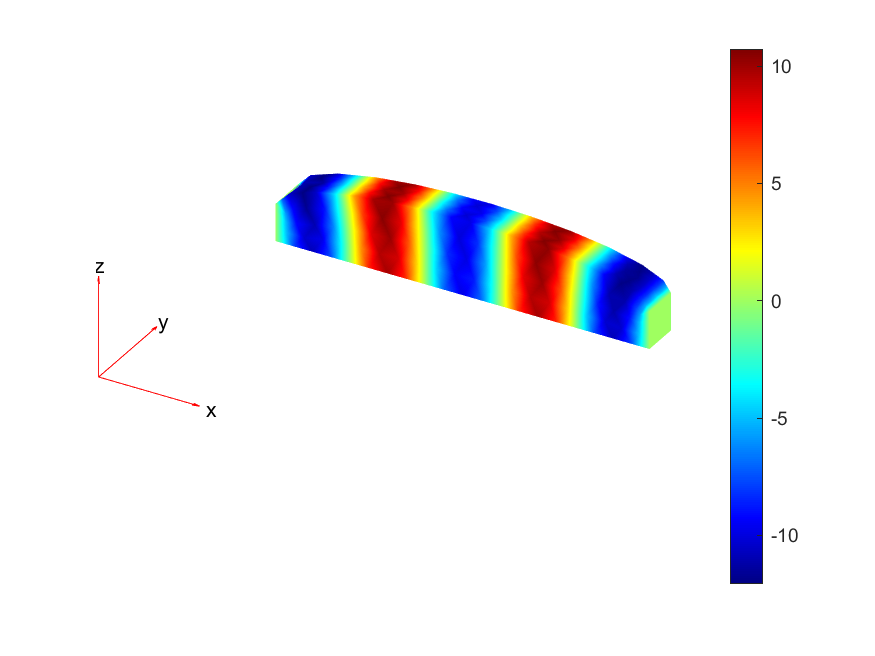}}
\scalebox{0.4}{\includegraphics{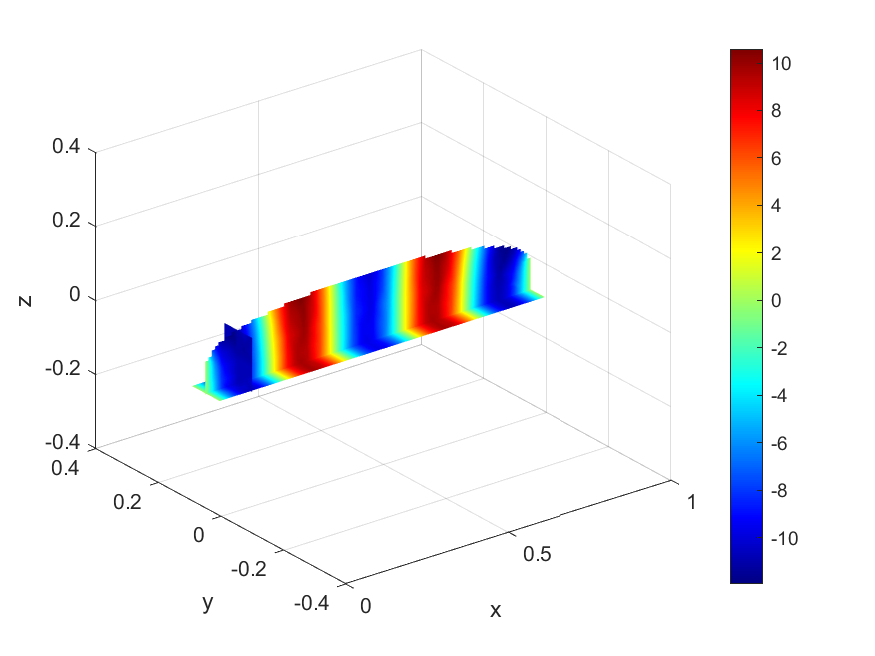}}
\scalebox{0.4}{\includegraphics{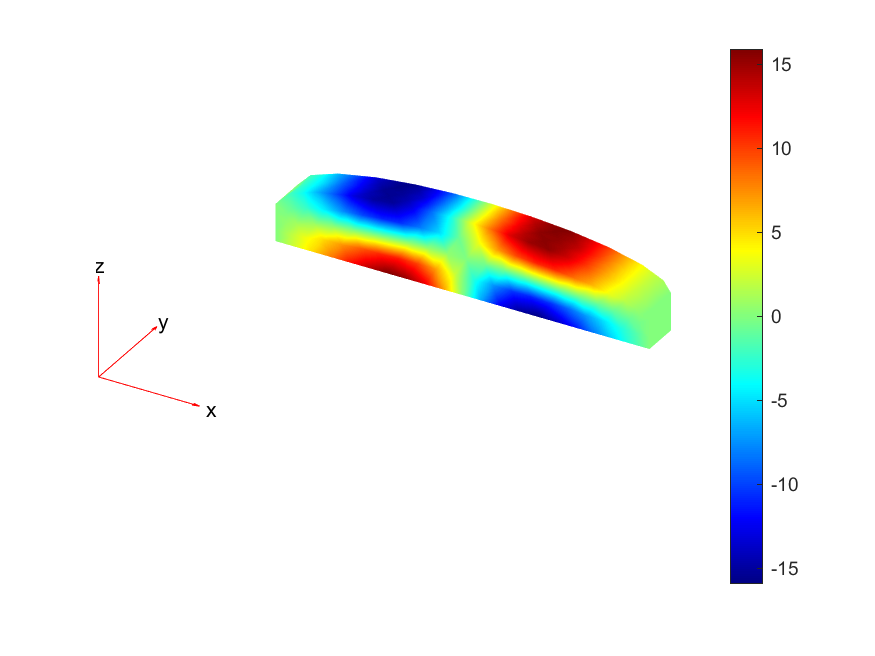}}
\scalebox{0.4}{\includegraphics{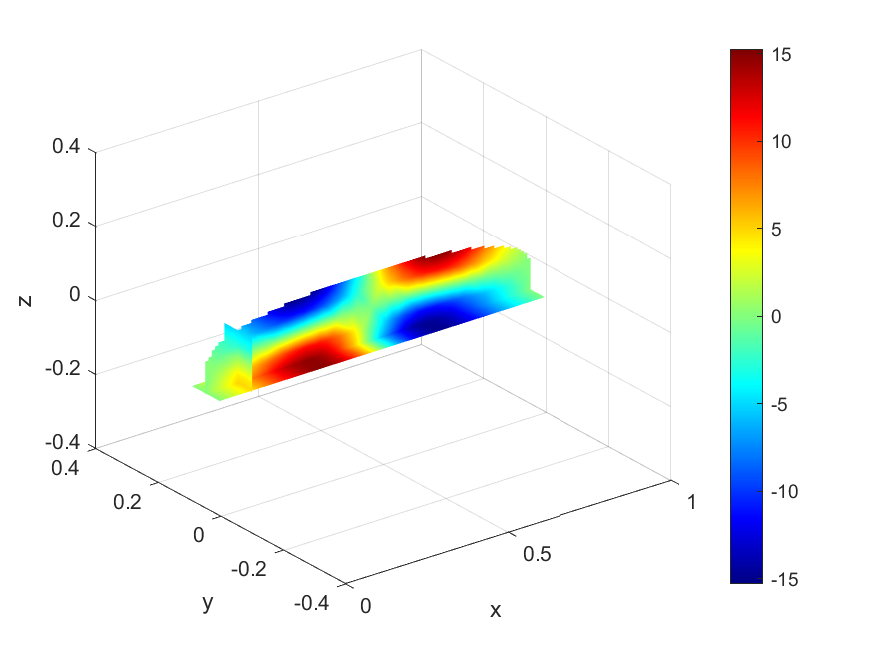}}
\end{center}
\caption{\small Approximations of eigenfunctions of \eqref{eq:ben4} with $G_\ee=(0,1)\times (0,\ee)\times (-\ee,\ee h(x_1))$ for a certain function $h$ satisfying \eqref{eq:ben13} and $\ee=0.1$. The figures are obtained choosing the eigenvalues $\lambda_1^\ee\approx 7.85,$
$\lambda_5^\ee\approx 236.43$ and  $\lambda_7^\ee\approx 326.39$. } \label{fig:ben4}
\end{figure}

Figure \ref{fig:ben4} shows numerical approximations of the eigenfunctions corresponding to the first,
fifth and   seventh  eigenvalue of \eqref{eq:ben4} when the domain  $G_\ee\!=\!(0,1)\times (0,\ee)\times (-\ee,\ee h(x_1))$ for a certain function $h$ satisfying \eqref{eq:ben13} and $\ee=0.1$ (see \eqref{eq:ben15}).
Here, $\lambda_1^\ee\approx 7.85,$ $\lambda_5^\ee\approx 236.43$ and  $\lambda_7^\ee\approx 326.39$. Note that the eigenfunction corresponding to the seventh eigenvalue presents transversal oscillations.

\section{The limit problem}\label{sec:ben4}
Because of the normalization condition \eqref{eq:ben6}  and the equation \eqref{eq:ben5}, for each fixed $n$ and  $\lambda^\ee=\lambda_n^\ee$ and $u^\ee=u_n^\ee$ the corresponding eigenfunction, we can write:
$$
 \int_{G_\ee} \vert\nabla u^\ee\vert^2  \, dx=\lambda^\ee \ee^2.$$ Introducing the change \eqref{eq:ben7}, we have
 $$
 \int_{G } \Big(\vert\partial_{y_1} u^\ee\vert^2  + \ee^{-2}\vert\partial_{y_2} u^\ee\vert^2 + \ee^{-2}\vert\partial_{y_3} u^\ee\vert^2 \Big) \,dy =\lambda^\ee \leq C_n, $$ where we have also used
  \eqref{eq:ben8}.   Thus, we get the bound
 \begin{equation}\label{eq:ben21} \Vert U^\ee\Vert_{H^1(G)}\leq C, \quad \Vert \partial_{y_2} U^\ee\Vert_{H^1(G)}\leq C\ee , \quad   \Vert\partial_{y_3}  U^\ee\Vert_{H^1(G)}\leq C \ee;  \end{equation}
 we note that, for a easy reading we have denoted $U^\ee(y):=u^\ee(y)$, and also the constants depend on the eigenvalue number $n$. Hence,
 for any sequence of $(\lambda^\ee, u^\ee)$, we can extract a subsequence, still denoted by $\ee$, such that
 \begin{equation}\label{eq:ben22}  \lambda^\ee\to \lambda^0, \quad  U^\ee  \to U^0 \mbox{ in } H^1(G)-weak, \quad  \mbox{ as  } \ee\to 0,\end{equation}
 for a certain  $\lambda^0>0$ and $U^0\in  H^1(G)$. Since, considering \eqref{eq:ben21},  the convergence
 $$ \partial_{y_i}   U^\ee  \to 0 \mbox{ in } L^2(G),  \mbox{ as  } \ee\to 0,  \mbox{ while } i=2,3,$$
also holds,  we deduce that $U^0(y_1,y_2,y_3)=U^0(y_1)$.

Hence, it remains to identify the limit problem satisfied by $(\lambda^0,U^0)$. We do it, by taking limits in \eqref{eq:ben5} for $v=\varphi(x_1)\in {\cal C}_0^\infty(\ell_0,\ell_1)$. Indeed, \eqref{eq:ben5} reads
$$
 \int_{G_\ee}  \partial_{x_1} u^\ee \varphi'   \, dx=\lambda^\ee  \int_{G_\ee} u^\ee \varphi \,dx,
$$
while introducing the  stretched coordinates, we have
$$
 \int_{G}  \partial_{y_1} U^\ee \varphi'   \, dy=\lambda^\ee  \int_{G} U^\ee \varphi \,dy.
$$

Now, according to \eqref{eq:ben22},  we pass to the   limits, as $\ee\to 0$, and get
\begin{equation}\label{eq:ben23}
 \int_{G }  \partial_{y_1} U^0 \varphi'   \, dy=\lambda^0  \int_{G } U^0 \varphi \,dy,\quad  \forall \varphi \in {\cal C}_0^\infty(\ell_0,\ell_1).
\end{equation}
Taking into account   \eqref{eq:ben1}, we rewrite \eqref{eq:ben23} as follows
$$
  \int_{\ell_0}^{\ell_1} \int_{D_{y_1}}  \partial_{y_1} U^0 \varphi'   \, dy_1 dy_2dy_3=\lambda^0   \int_{\ell_0}^{\ell_1} \int_{D_{y_1}}  U^0 \varphi \,dy_1 dy_2dy_3,\quad  \forall \varphi \in {\cal C}_0^\infty(\ell_0,\ell_1).
$$
Now, using a density argument, we  obtain
\begin{equation*}
  \int_{\ell_0}^{\ell_1} \vert D_{y_1}\vert \,  \partial_{y_1} U^0 \varphi'   \, dy_1 =\lambda^0   \int_{\ell_0}^{\ell_1} \vert D_{y_1} \vert  \, U^0 \varphi \,dy_1 ,\quad  \forall \varphi \in H_0^1(\ell_0,\ell_1),
\end{equation*}
namely, the weak formulation of the following Dirichlet problem
\begin{equation}\label{eq:ben24}
\left\{
\begin{array}{ll}-\partial_{x_1}\Big( \vert D_{x_1} \vert \, \partial_{x_1} U^0\Big) =\lambda^0 \vert D_{x_1} \vert  \,
U^0, \quad x_1\in (\ell_0,\ell_1), \vspace{0.2cm}\\  U^0(\ell_0)=0,\quad U^0(\ell_1)=0.
\end{array}\right.
\end{equation}
which takes into account the geometry of $G_\ee$.

 In addition to the above proofs, we note that the convergence \eqref{eq:ben22}, also implies the convergence of the eigenfunctions in $L^2(G)$, and consequently \eqref{eq:ben6} ensures $U^0\not\equiv0$.

Hence we have proved the following result:
\begin{theorem}\label{theorem:ben1}
Assume that the area of the transverse sections of the domain $G$, $\vert D_{x_1}\vert$ is
a stepwise continuous function for $x_1\in (\ell_0,\ell_1)$ satisfying \eqref{eq:ben2}. Then, for any fixed $n\in \mathbb{N}$, and for any subsequence of $\ee$ still denoted by $\ee$, such that \begin{equation*}  \lambda_n^\ee\to \lambda_n^0, \quad  U_n^\ee  \to U_n^0 \quad  \mbox{ in } H^1(G)-weak, \quad  \mbox{ as  } \ee\to 0,\end{equation*}
we have that  $\lambda_n^0$ is an eigenvalue of \eqref{eq:ben24} and $U_n^0$ is an associated eigenfunction.\end{theorem}

 Let us analyze the result of Theorem \ref{theorem:ben1} for different geometries of the domain $G_\ee$ described in Section \ref{sec:ben2_1}.

In the case where the domain $G_\ee$ is given by \eqref{eq:ben9} the computations of eigenvalues and eigenfunctions can be done explicitly and they are in Section \ref{sec:ben3} (cf. Figure \ref{fig:ben2}).

In the case where $G_\ee$ is given by  \eqref{eq:ben11} with $\ell_1>0$, the function $\vert D_{x_1}\vert $ is stepwise continuous function defined by
$$\vert D_{x_1}\vert =4 \, \mbox{ when } x_1\in [0,2^{-1}\ell_1]\quad  \mbox{ and } \quad  \vert D_{x_1}\vert =1  \, \mbox{ when } x_1\in  (2^{-1}\ell_1,\ell_1].$$
Let us refer to \cite{PanasenkoPerez2007} to compare with  a 2D domain and a different technique of approach, and also notice that the above proofs  (and hence the  result of Theorem \ref{theorem:ben1})  apply to less smooth functions $\vert D_{x_1}\vert$.

Finally, for simplicity we consider the case when  $G_\ee$ is given by \eqref{eq:ben12}, where  $\ell_0=0$ and the function   $\vert D_{x_1}\vert $ reads
 $$\vert D_{x_1}\vert =h(x_1), \,\,  x_1\in (0,\ell_1), $$
 with $h$ a smooth function satisfying \eqref{eq:ben13}.

To complete the convergence results  obtaining the convergence of the spectrum of  \eqref{eq:ben4} towards that of \eqref{eq:ben24} with conservation of the multiplicity we can use a convergence result  for  operators on Hilbert spaces, both depending on a small parameter  (cf.  Section~III.1 of \cite{OlShYo92}).
The application of such a result to the problem under consideration is left as an open problem to be addressed  by the authors in forthcoming publication.

\section{Appendix: different boundary conditions}\label{sec:ben5}

In this section we  present illustrative computations for the case of the spectrum of  Neumann Laplacian, namely, of problem
\begin{equation}\label{eq:ben25}
\left\{\begin{array}{ll}-\Delta u^\ee=\lambda^\ee u^\ee \hbox{ in }G_\ee, \vspace{0.2cm}\\
\displaystyle{ \frac{\partial u^\ee}{\partial\nu}=0 \hbox{ on } \partial G_\ee,}
\end{array}\right.
\end{equation}
and the spectrum of Dirichlet Laplacian, namely, problem
\begin{equation}\label{eq:ben26}
\left\{\begin{array}{ll}-\Delta u^\ee=\lambda^\ee u^\ee \hbox{ in }G_\ee, \vspace{0.2cm}\\
u^\ee=0 \hbox{ on } \partial G_\ee.
\end{array}\right.
\end{equation}

For the Neumann problem \eqref{eq:ben25} posed in the prism $G_\ee$ defined by \eqref{eq:ben9}, we can compute explicitly the eigenvalues and the corresponding eigenfunctions which are given by
\begin{equation}\label{eq:ben27}
\lambda_{mrs}^\ee =\Big(\dfrac{m\pi}{\ell_1}\Big)^2+\Big(\dfrac{r\pi}{\ee}\Big)^2 + \Big(\dfrac{s\pi}{\ee}\Big)^2, \quad  m, r, s \in \mathbb{N}\cup\{0\},
\end{equation}
\begin{equation}\label{eq:ben28}
u_{mrs}^\ee =A_{mrs}\cos\Big(\dfrac{m\pi x_1}{\ell_1}\Big)\cos\Big(\dfrac{r\pi x_2}{\ee}\Big)\cos\Big(\dfrac{s\pi x_3}{\ee}\Big),  \quad A_{mrs}\in\mathbb{R},  \,\,   m, r, s \in \mathbb{N}\cup\{0\}.
\end{equation}
Now, the separation of variables leads to the problems \begin{eqnarray*}
&F''(x_1)+\mu_1 F(x_1)=0  \quad x_1 \in (0,\ell_1),\quad  F'(0)=F'(\ell_1)=0, \vspace{0.25cm} \end{eqnarray*}
\eqref{eq:ben17} and \eqref{eq:ben18}.
Here, the first eigenvalue is equal to zero,  $\lambda_0^\ee=0$, and the corresponding eigenfunctions are the constants.

\begin{figure}[t]
\begin{center}
\scalebox{0.4}{\includegraphics{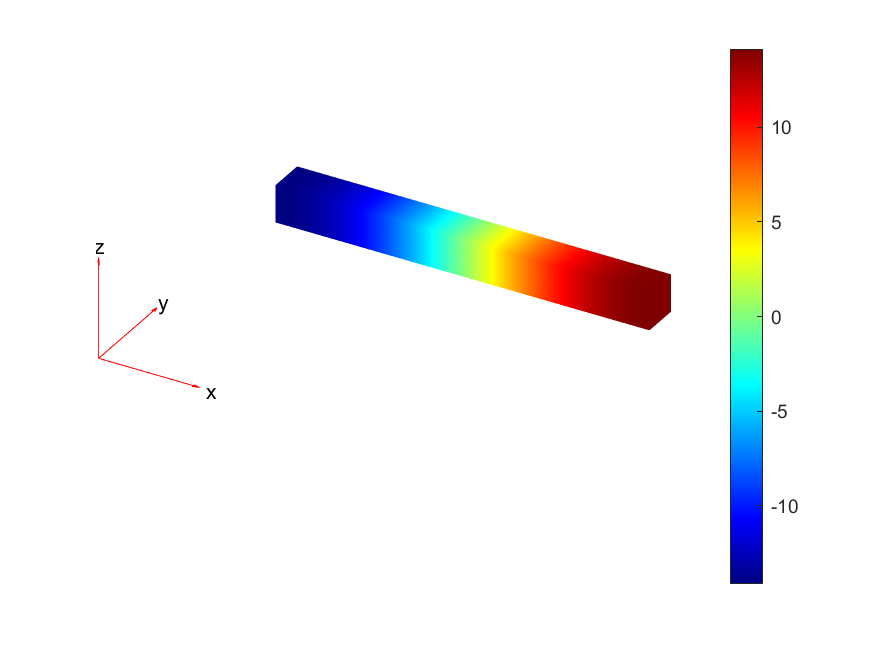}}
\scalebox{0.4}{\includegraphics{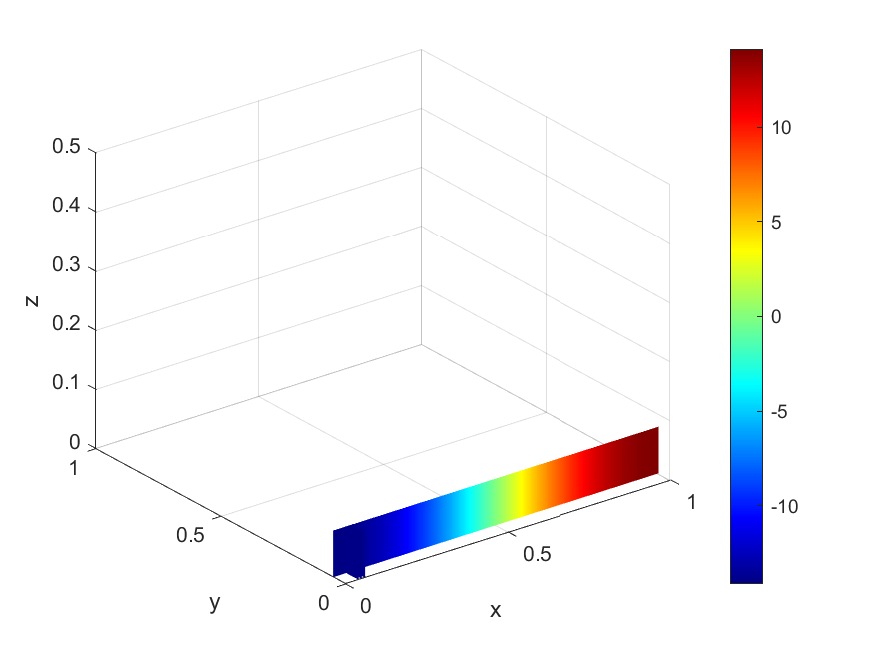}}
\scalebox{0.4}{\includegraphics{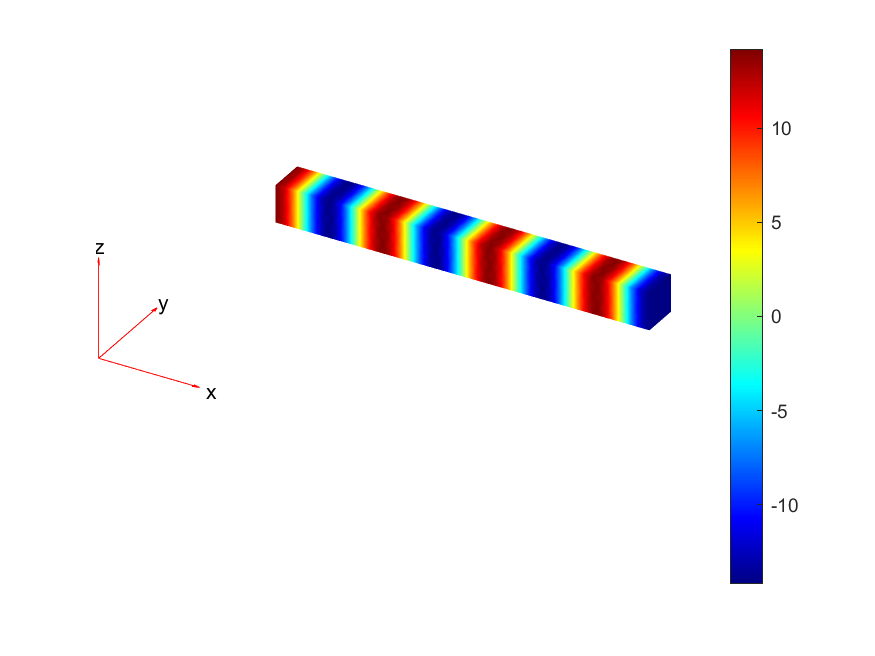}}
\scalebox{0.4}{\includegraphics{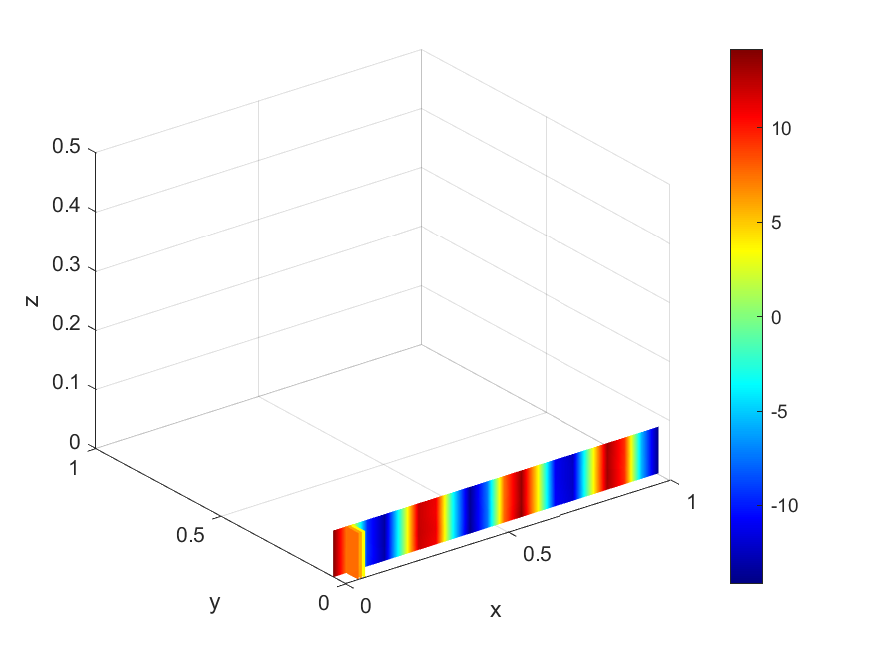}}
\scalebox{0.4}{\includegraphics{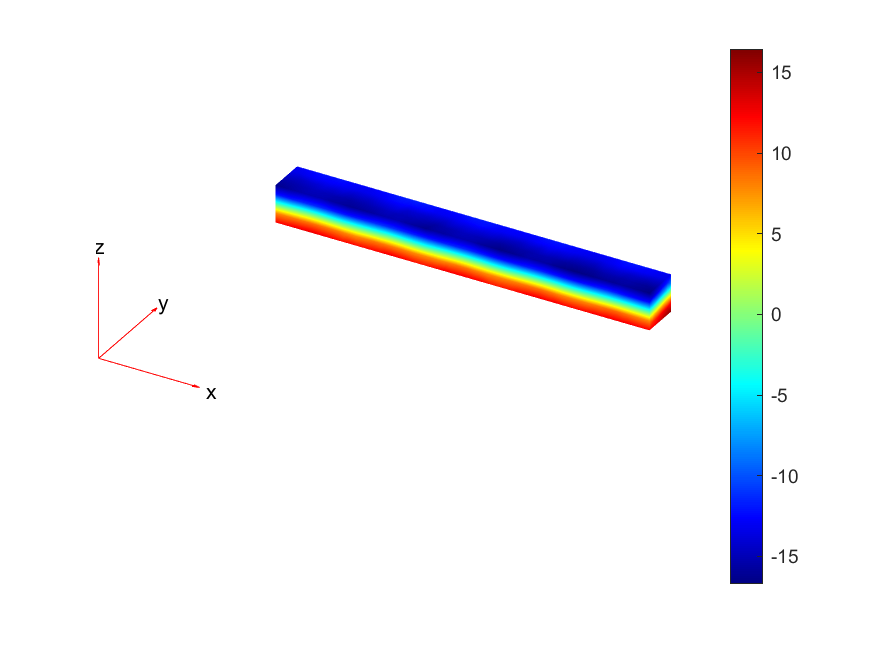}}
\scalebox{0.4}{\includegraphics{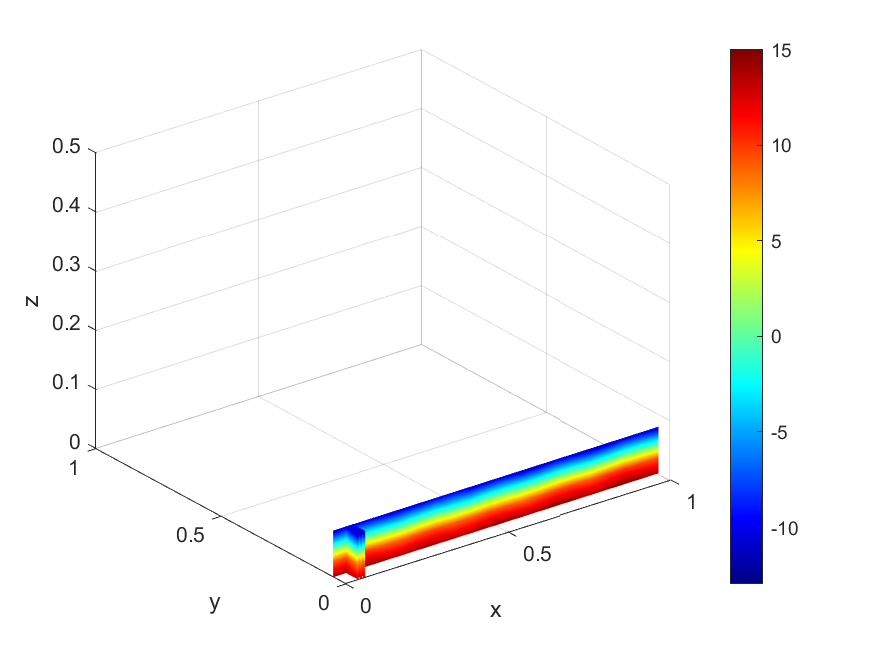}}
\end{center}
\caption{\small  Approximations of eigenfunctions of \eqref{eq:ben25} with $G_\ee=(0,1)\times (0,\ee)\times (0,\ee)$ and $\ee=0.1$. The figures are obtained choosing the eigenvalues $\lambda_1^\ee=\pi^2\approx 9.87,$ $\lambda_7^\ee=49\pi^2\approx 484.09$ and $\l_{11}^\ee=100\pi^2\approx 992.65.$ } \label{fig:ben5}
\end{figure}

Similarly to problem \eqref{eq:ben4}, from \eqref{eq:ben27} and \eqref{eq:ben28}, we observe that, for $\ee$ small, the oscillations of the eigenfunctions corresponding to the low frequencies
are longitudinal (associated again with the parameters $r=s=0$) and in order to capture transverse  oscillations of the eigenfunctions we need to deal with the high frequencies (cf. Figure \ref{fig:ben5}).

Figure \ref{fig:ben5} shows numerical approximations of the eigenfunctions corresponding to the first, seventh and eleventh eigenvalue of \eqref{eq:ben25} different from zero when the domain $G_\ee$ is a prism $G_\ee=(0,1)\times (0,\ee)\times (0,\ee)$ and $\ee=0.1$ (see \eqref{eq:ben9}). Here, $\lambda_1^\ee=\pi^2\approx 9.87,$ $\lambda_7^\ee=49\pi^2\approx 484.08$ and it is necessary to reach the eleventh eigenvalue $\l_{11}^\ee=100\pi^2\approx 992.65$ to capture transversal oscillations.

\begin{figure}[t]
\begin{center}
\scalebox{0.4}{\includegraphics{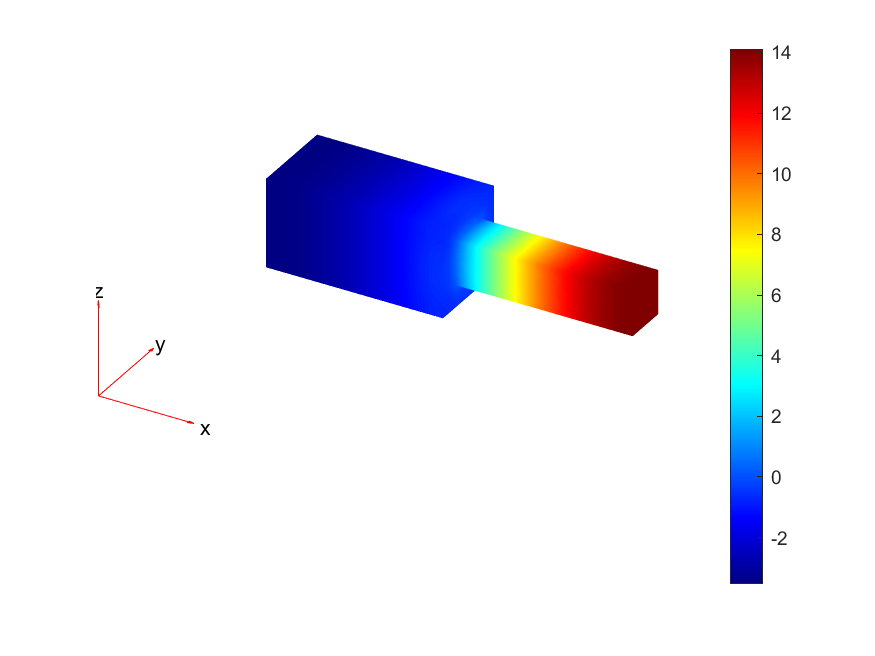}}
\scalebox{0.4}{\includegraphics{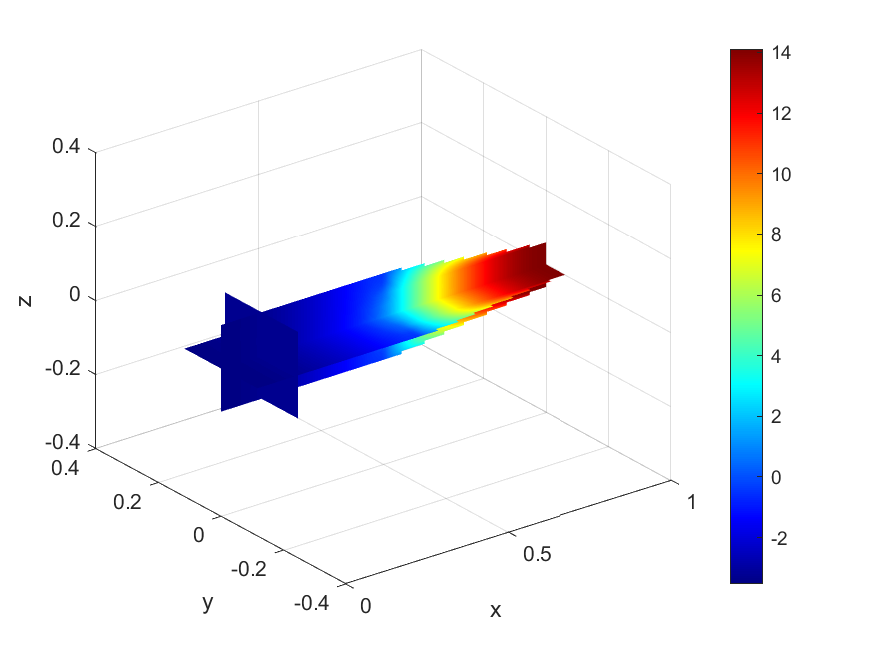}}
\scalebox{0.4}{\includegraphics{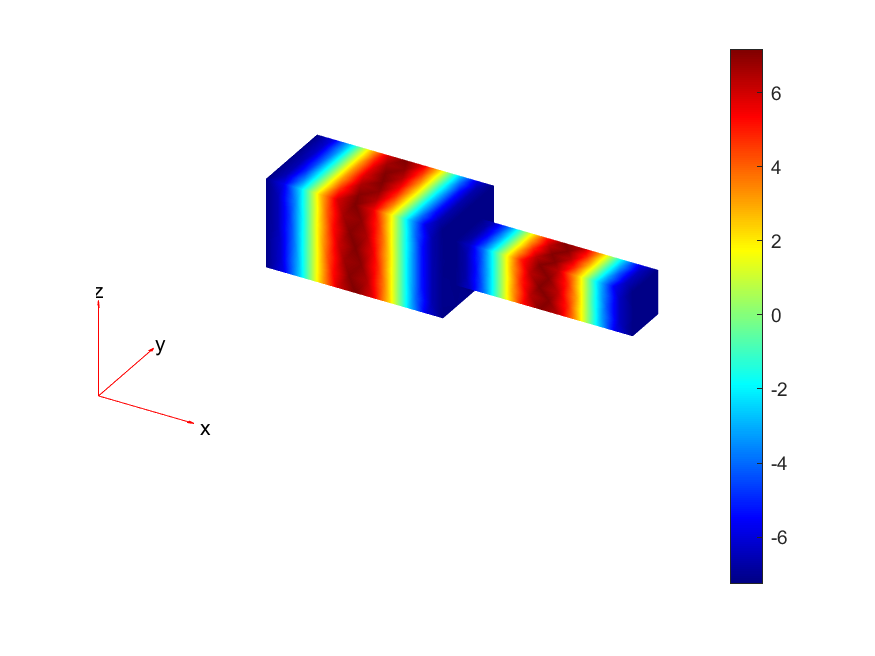}}
\scalebox{0.4}{\includegraphics{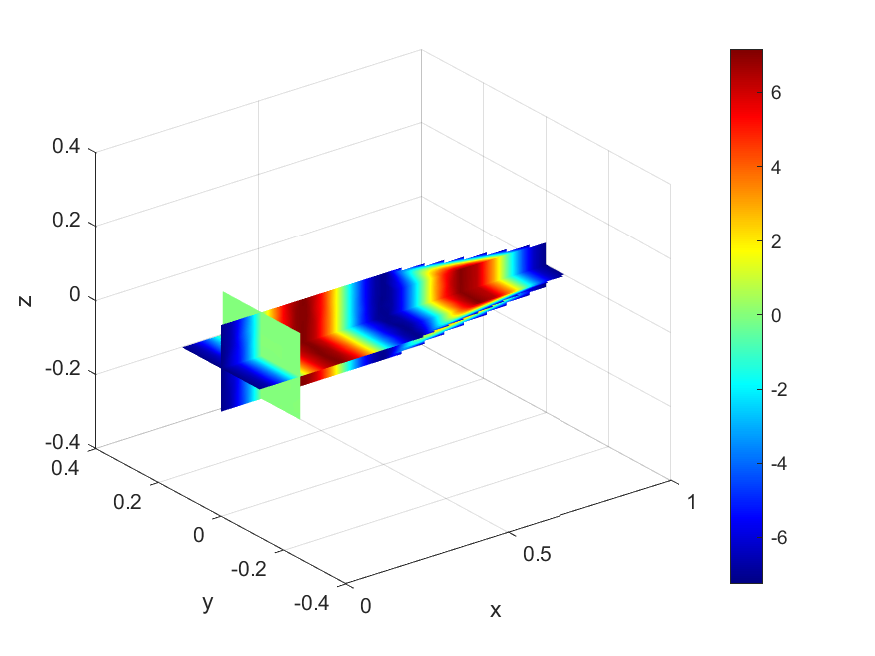}}
\scalebox{0.4}{\includegraphics{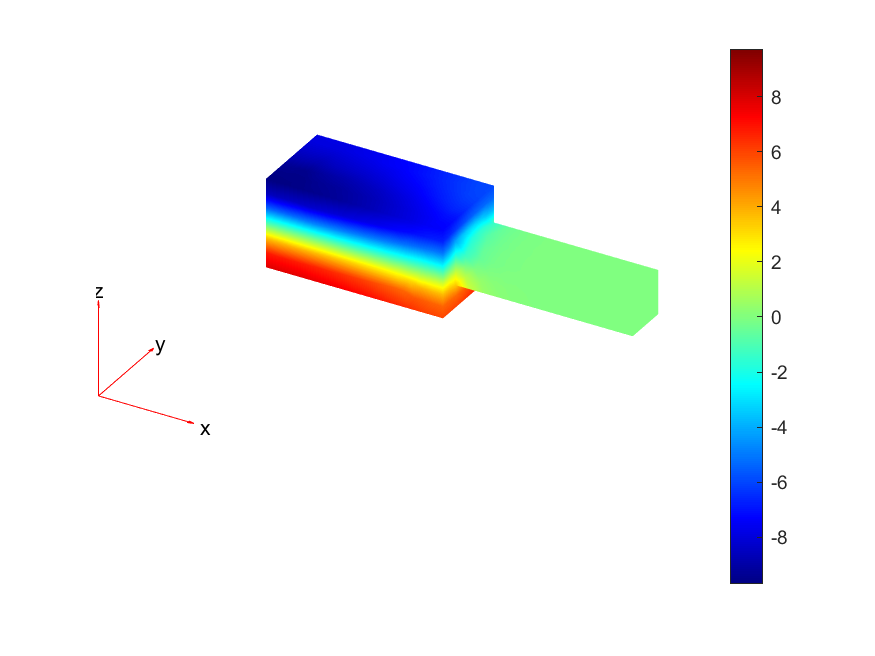}}
\scalebox{0.4}{\includegraphics{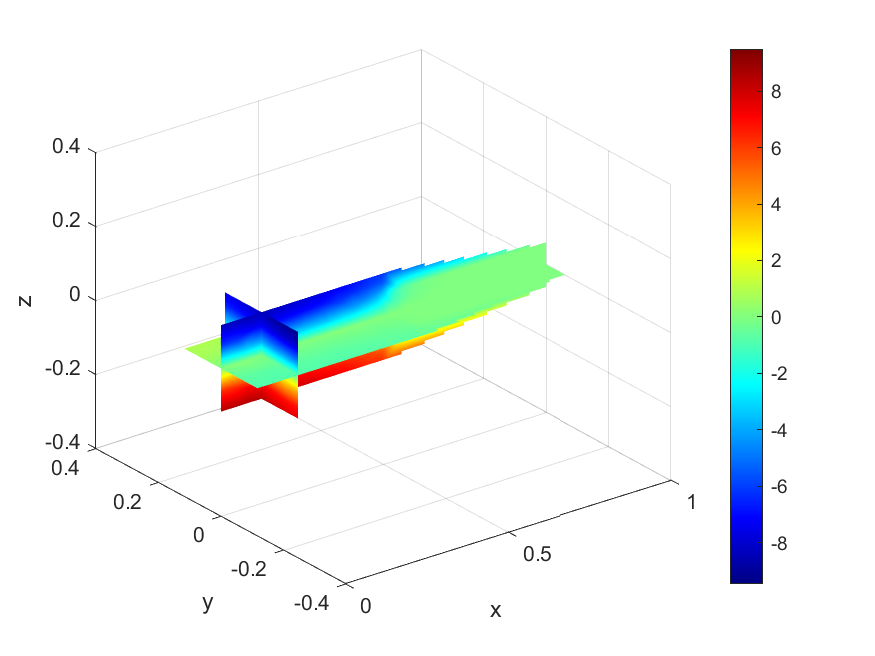}}
\end{center}
\caption{\small Approximations of eigenfunctions of \eqref{eq:ben25} with $G_\ee=(0,2^{-1})\times(-\ee,\ee)\times(-\ee,\ee) \, \cup \,
\{2^{-1}\}\times (-\ee,\ee)\times(-\ee,\ee) \, \cup \, (1/2,1)\times (-\ee2^{-1},\ee2^{-1})\times (-\ee2^{-1},\ee2^{-1})$ and $\ee=8^{-1}$. The figures are obtained choosing the eigenvalues $\lambda_1^\ee \approx 9.35,$ $\lambda_4^\ee\approx 158.11$, $\lambda_5^\ee\approx 160.85$.}\label{fig:ben6}
\end{figure}

Figure \ref{fig:ben6} shows numerical  approximations of the eigenfunctions corresponding to the first, fourth and  fifth eigenvalue of \eqref{eq:ben25} different from zero when the domain $G_\ee=(0,2^{-1})\times(-\ee,\ee)\times(-\ee,\ee) \, \cup \,
\{2^{-1}\}\times (-\ee 2^{-1},\ee 2^{-1})\times(-\ee2^{-1},\ee2^{-1}) \, \cup \, (2^{-1},1)\times (-\ee2^{-1},\ee2^{-1})\times (-\ee2^{-1},\ee2^{-1})$ and $\ee=8^{-1}$ (see \eqref{eq:ben11}). Here, $\lambda_1^\ee \approx 9.35,$ $\lambda_4^\ee\approx 158.11$, $\lambda_5^\ee\approx 160.85$.
We already notice certain phenomena of localization for the eigenfunctions which need  a thorough study.

\begin{figure}[t]
\begin{center}
\scalebox{0.4}{\includegraphics{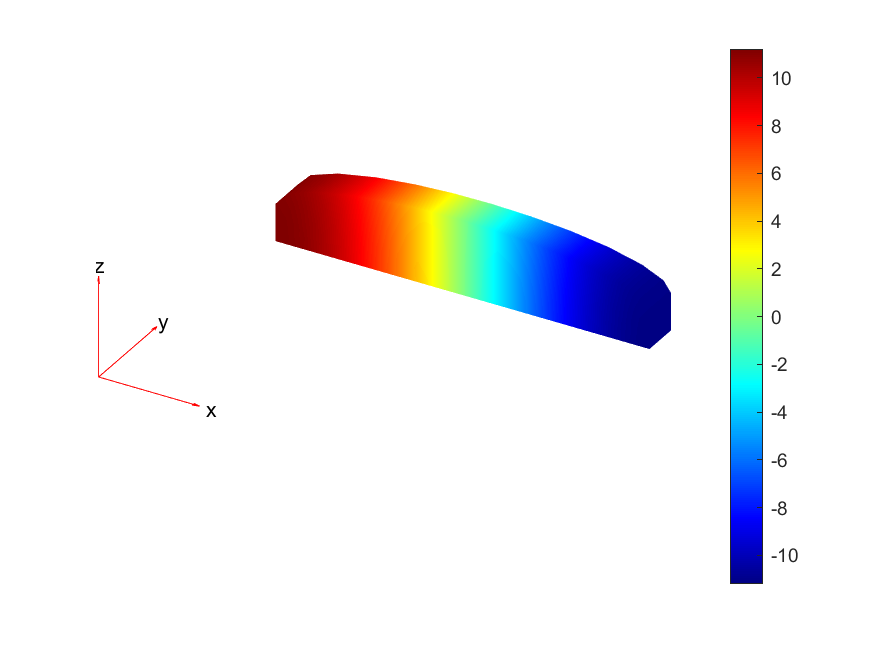}}
\scalebox{0.4}{\includegraphics{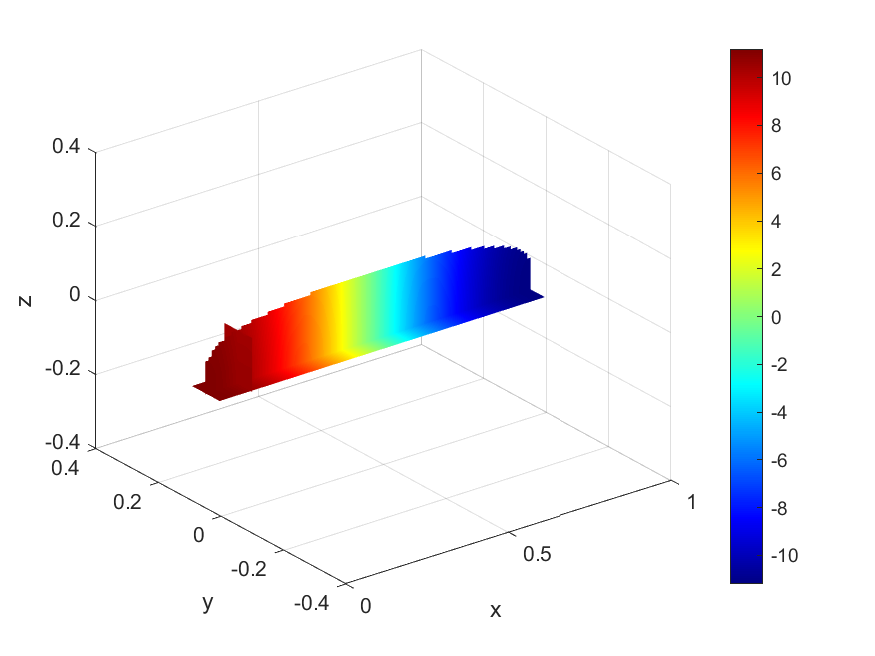}}
\scalebox{0.4}{\includegraphics{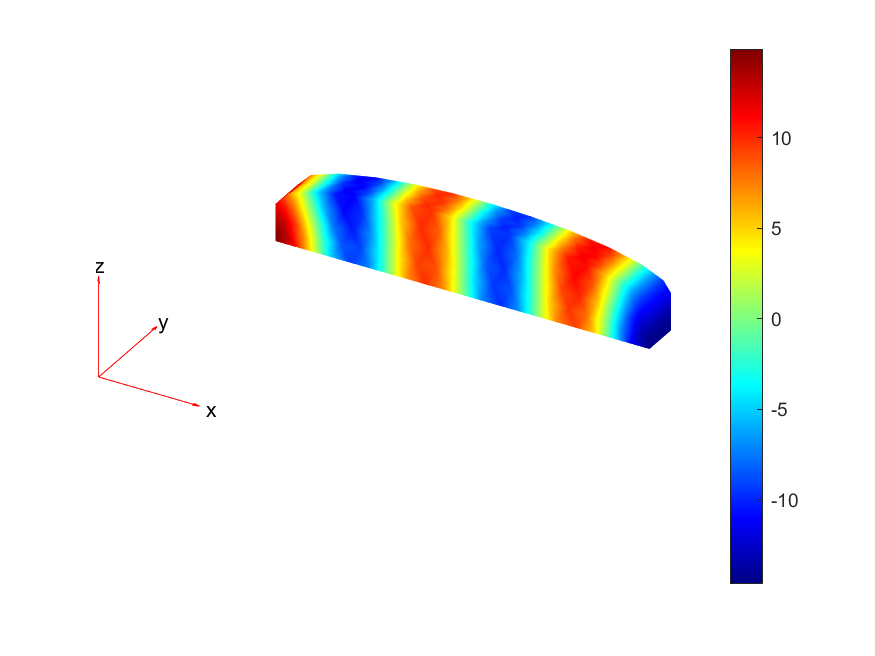}}
\scalebox{0.4}{\includegraphics{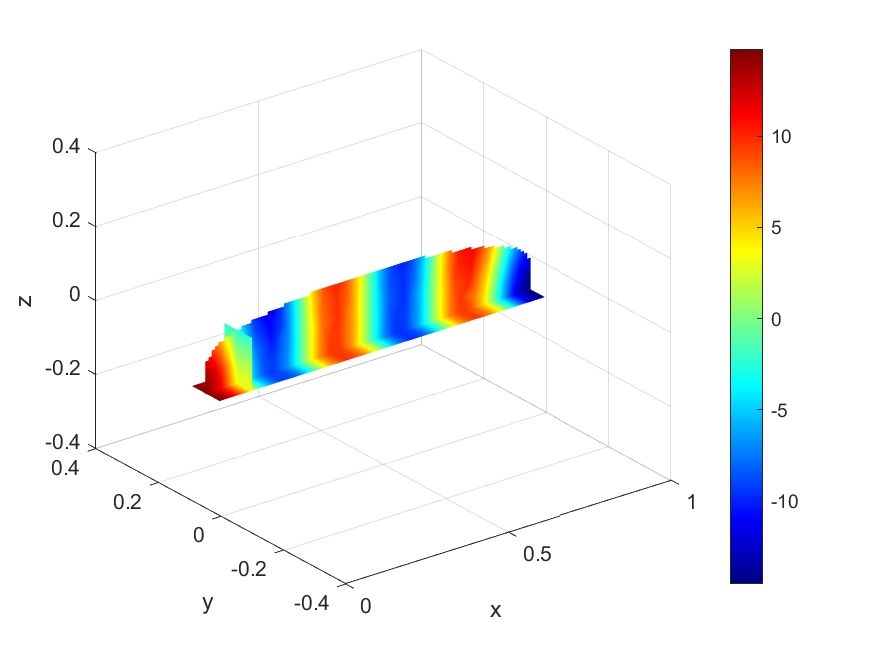}}
\scalebox{0.4}{\includegraphics{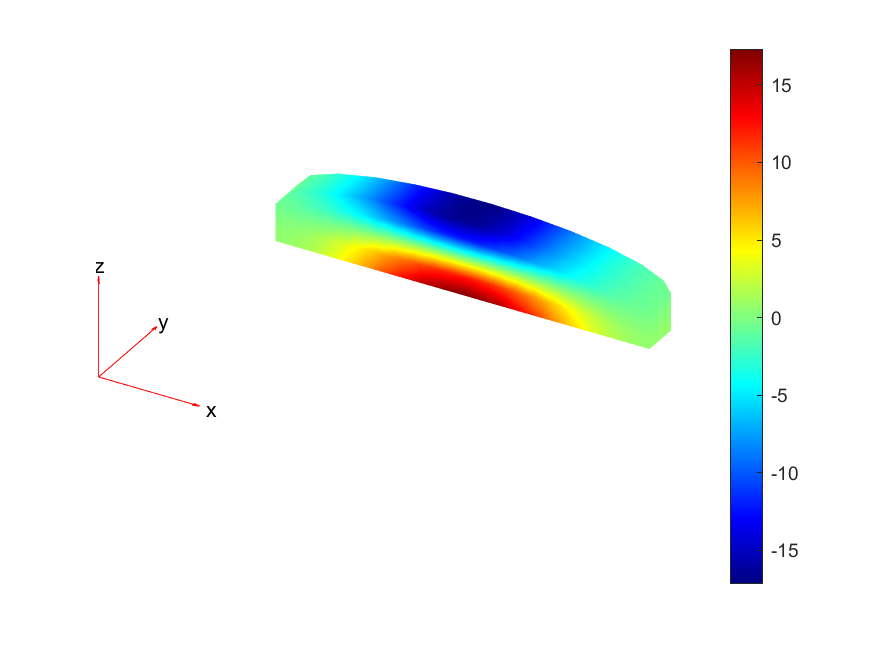}}
\scalebox{0.4}{\includegraphics{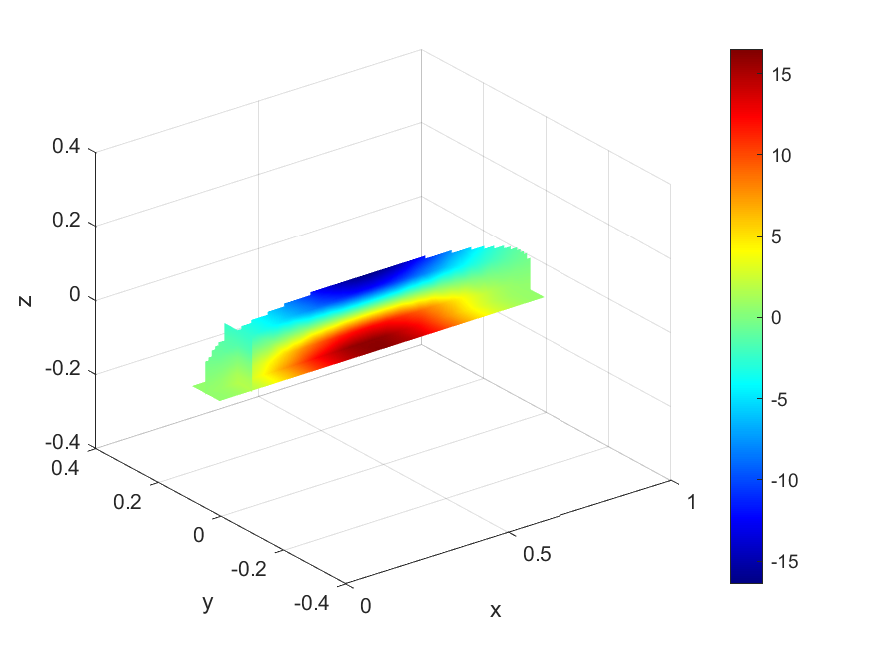}}
\end{center}
\caption{\small Approximations of eigenfunctions of \eqref{eq:ben25} with $G_\ee=(0,1)\times (0,\ee)\times (-\ee,\ee h(x_1))$ for a certain function $h$ satisfying \eqref{eq:ben13} and $\ee=0.1$. The figures are obtained choosing the eigenvalues $\lambda_1^\ee\approx 11.62,$ $\lambda_5^\ee\approx 249.49$ and  $\lambda_6^\ee\approx 271.67$. } \label{fig:ben7}
\end{figure}

Figure \ref{fig:ben7} shows numerical approximations of the eigenfunctions corresponding to the first, fifth and   sixth  eigenvalue of \eqref{eq:ben25} different from zero when the domain $G_\ee=(0,1)\times (0,\ee)\times (-\ee,\ee h(x_1))$ for a certain function $h$ satisfying \eqref{eq:ben13} and $\ee=0.1$ (see \eqref{eq:ben15}).
Here, $\lambda_1^\ee\approx 11.62,$ $\lambda_5^\ee\approx 249.49$ and  $\lambda_6^\ee\approx 271.67$.

For the Dirichlet problem \eqref{eq:ben26} posed in the prism $G_\ee$ defined by \eqref{eq:ben9}, by separation of variables, we can compute explicitly the eigenvalues and the corresponding eigenfunctions which are given by
\begin{equation}\label{eq:ben29}
\lambda_{mrs}^\ee =\Big(\dfrac{m\pi}{\ell_1}\Big)^2+\Big(\dfrac{r\pi}{\ee}\Big)^2 + \Big(\dfrac{s\pi}{\ee}\Big)^2, \quad  m, r, s \in \mathbb{N},
\end{equation}
\begin{equation}\label{eq:ben30}
u_{mrs}^\ee =A_{mrs}\sin\Big(\dfrac{m\pi x_1}{\ell_1}\Big)\sin\Big(\dfrac{r\pi x_2}{\ee}\Big)\sin\Big(\dfrac{s\pi x_3}{\ee}\Big),  \quad A_{mrs}\in\mathbb{R},  \,\,   m, r, s \in \mathbb{N}.
\end{equation}
Now, the separation of variables leads to the problems \eqref{eq:ben16},
\begin{eqnarray*}
&G''(x_2)+\mu_2^\ee F(x_2)=0  \quad x_2 \in (0,\ee),\quad  G(0)=G(\ee)=0, \vspace{0.25cm} \\
&H''(x_3)+\mu_3^\ee H(x_3)=0  \quad x_3 \in (0,\ee),\quad  H(0)=H(\ee)=0.
\end{eqnarray*}

\begin{figure}[t]
\begin{center}
\scalebox{0.4}{\includegraphics{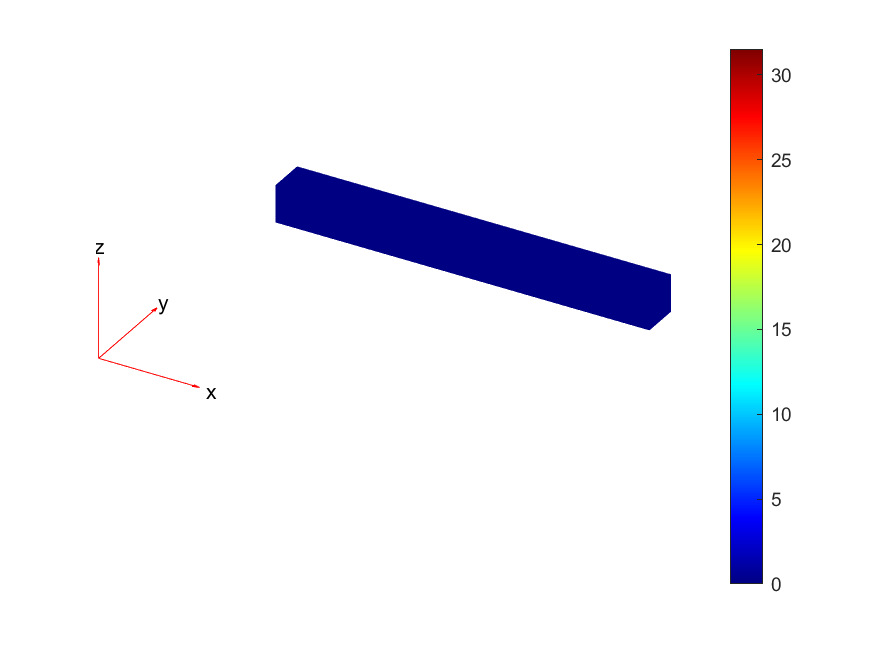}}
\scalebox{0.4}{\includegraphics{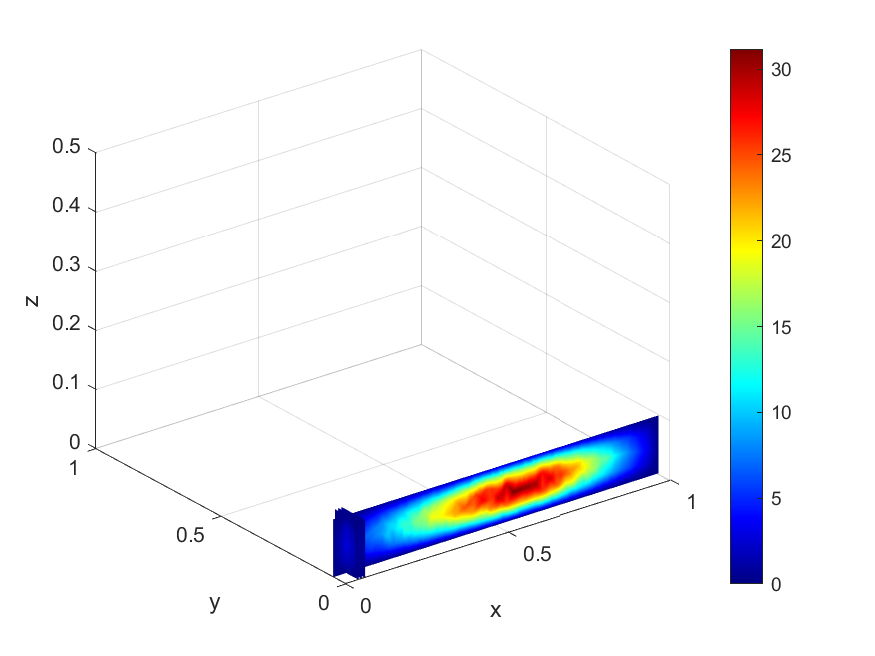}}
\scalebox{0.4}{\includegraphics{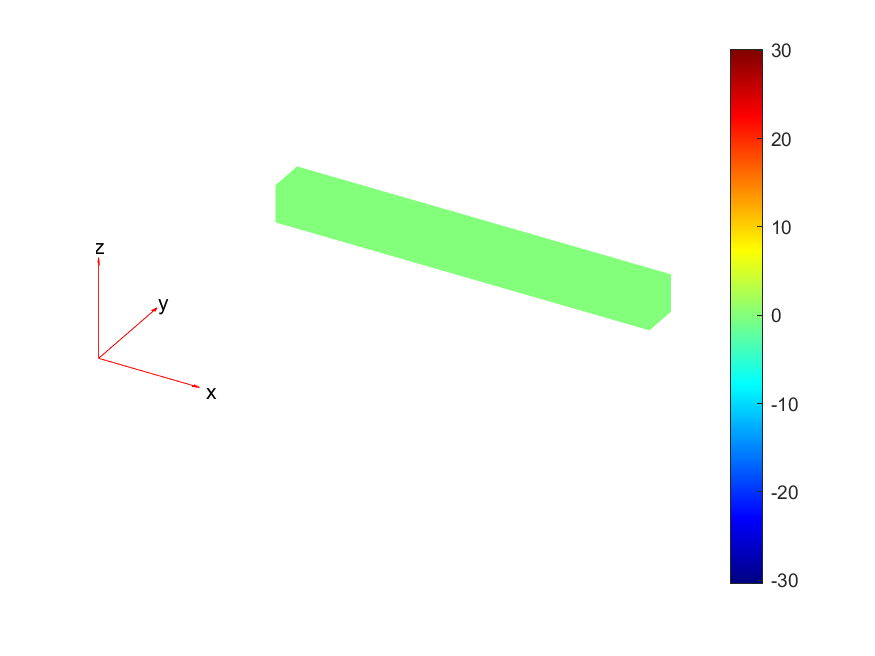}}
\scalebox{0.4}{\includegraphics{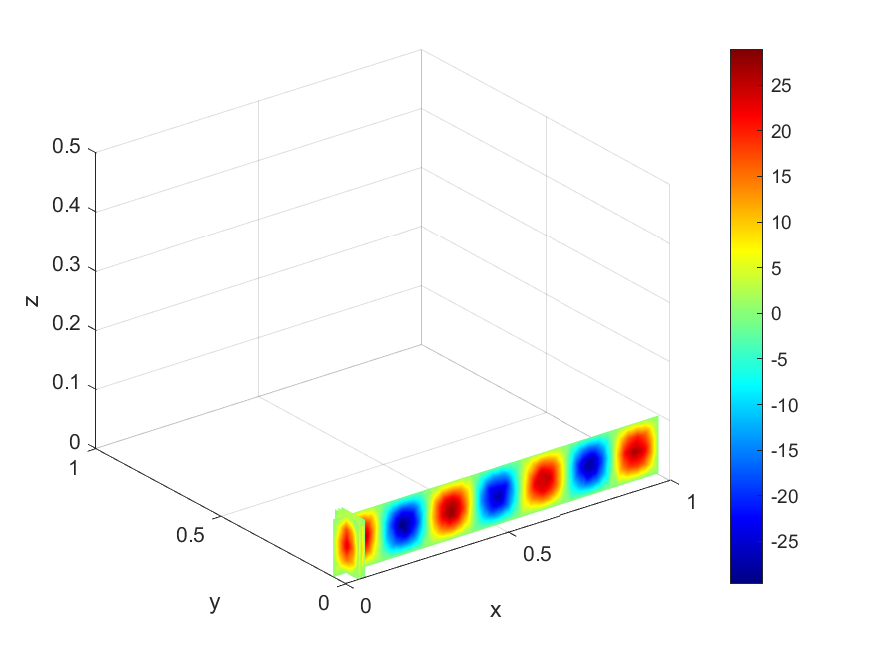}}
\scalebox{0.4}{\includegraphics{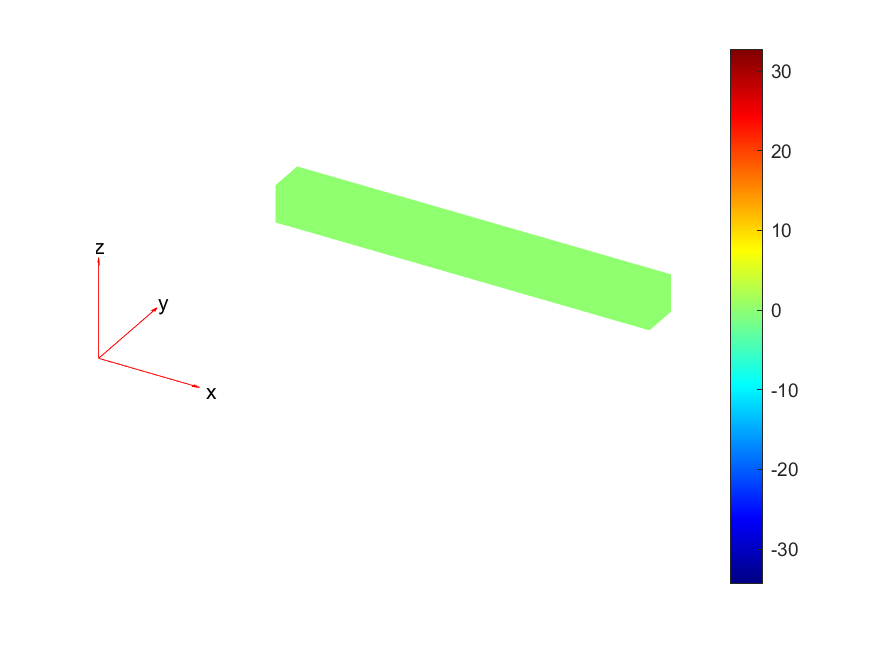}}
\scalebox{0.4}{\includegraphics{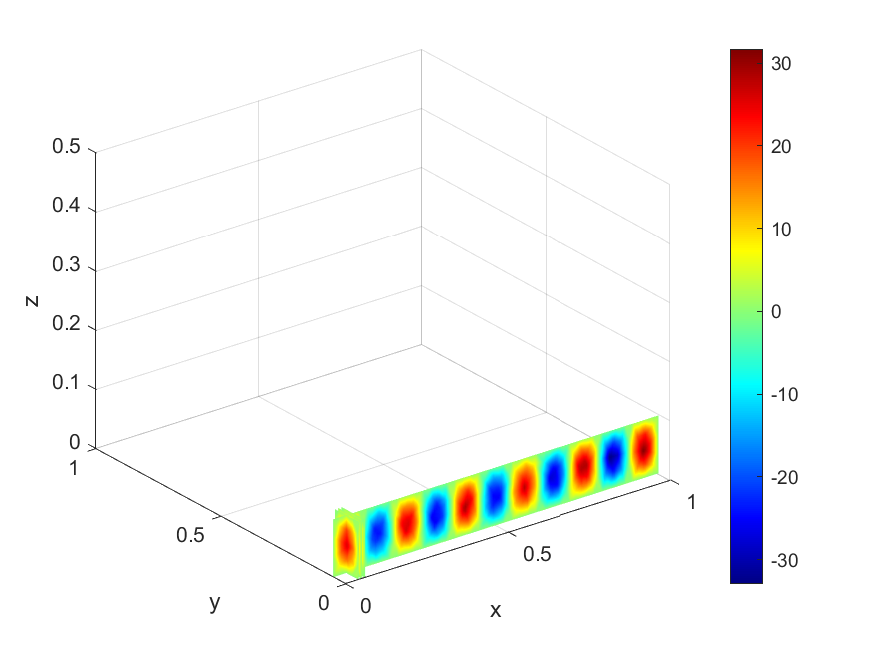}}
\end{center}
\caption{\small Approximations of eigenfunctions of \eqref{eq:ben26} with $G_\ee=(0,1)\times (0,\ee)\times (0,\ee)$ and $\ee=0.1$. The figures are obtained choosing the eigenvalues $\lambda_1^\ee=201\pi^2\approx 2027.60,$ $\lambda_7^\ee=249\pi^2\approx 2540.44$ and $\l_{11}^\ee=321\pi^2\approx 3336.87.$ } \label{fig:ben8}
\end{figure}

Contrary to problems \eqref{eq:ben4} and \eqref{eq:ben25}, from \eqref{eq:ben29} and \eqref{eq:ben30}, we observe that now, for $\ee$ small, the oscillations of the eigenfunctions corresponding to the low frequencies
are transversal, now associated with the parameters $m, r, s$ different from zero;
cf. Figure \ref{fig:ben8}.

Figure \ref{fig:ben8} shows numerical approximations of the eigenfunctions corresponding to the first, seventh and eleventh eigenvalue of \eqref{eq:ben26} when the domain $G_\ee$ is a prism $G_\ee=(0,1)\times (0,\ee)\times (0,\ee)$ and $\ee=0.1$ (see \eqref{eq:ben9}). Here, $\lambda_1^\ee=201\pi^2\approx 2027.60,$ $\lambda_7^\ee=249\pi^2\approx 2540.44$ and $\l_{11}^\ee=321\pi^2\approx 3336.87.$

\begin{figure}[t]
\begin{center}
\scalebox{0.4}{\includegraphics{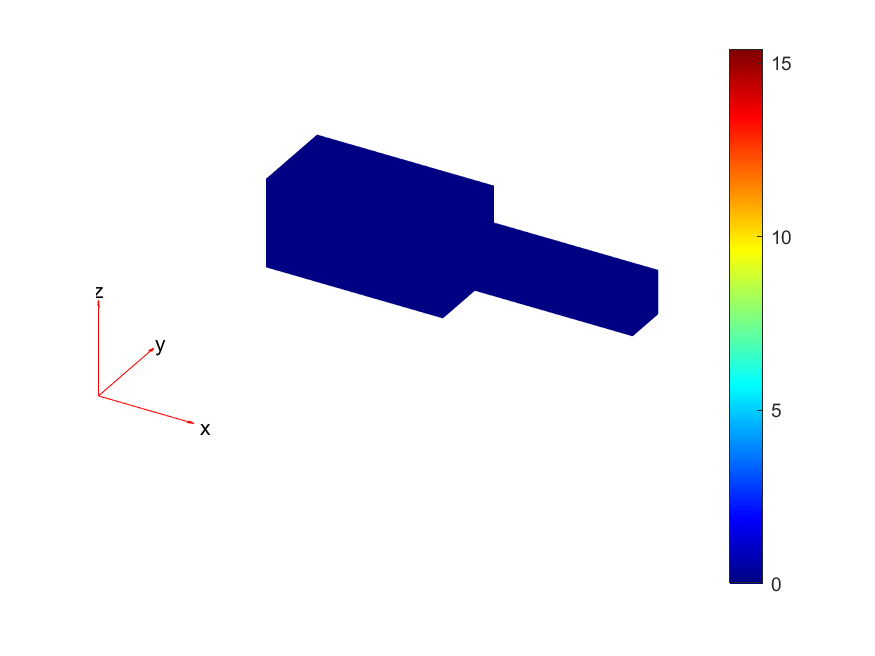}}
\scalebox{0.4}{\includegraphics{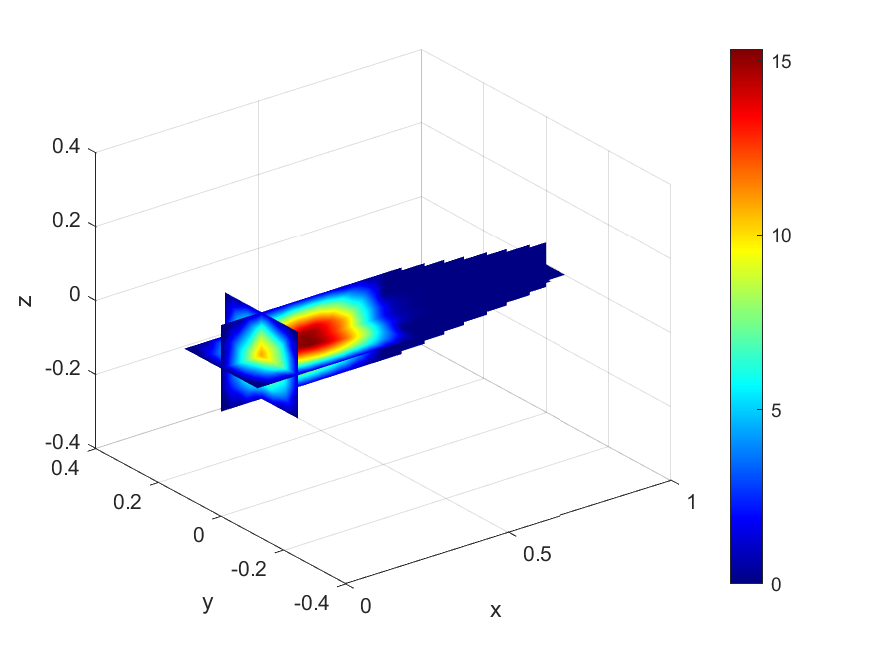}}
\scalebox{0.4}{\includegraphics{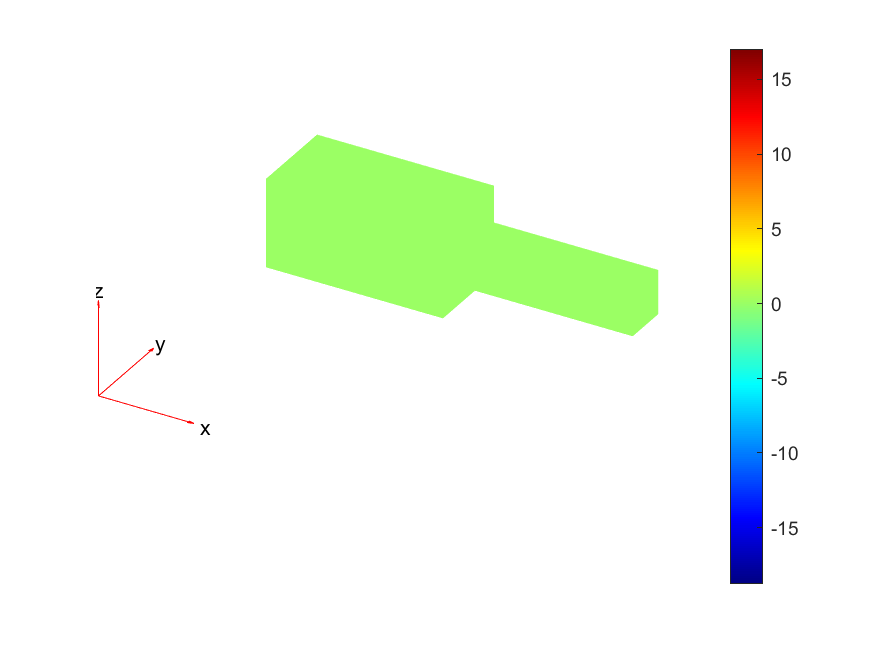}}
\scalebox{0.4}{\includegraphics{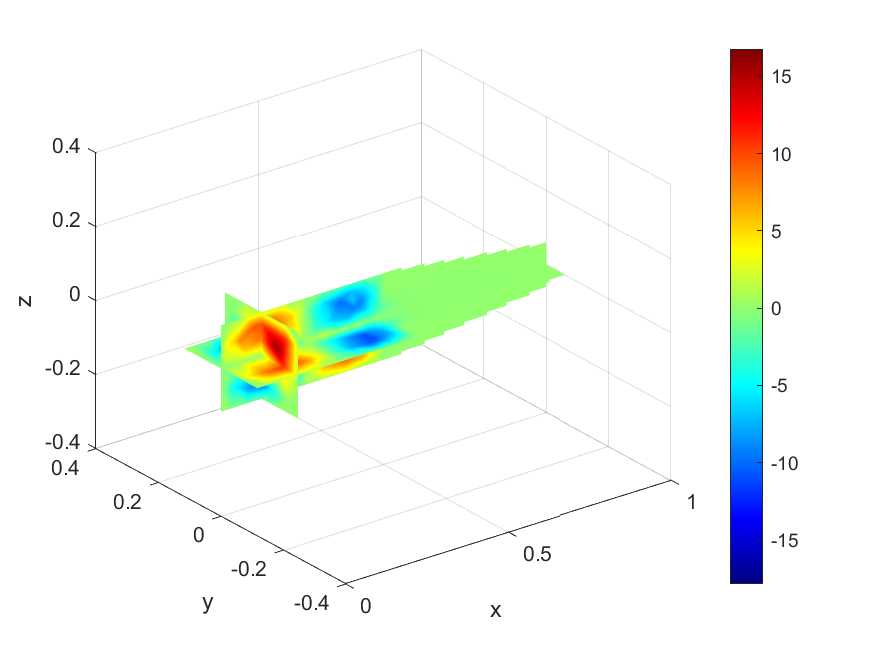}}
\scalebox{0.4}{\includegraphics{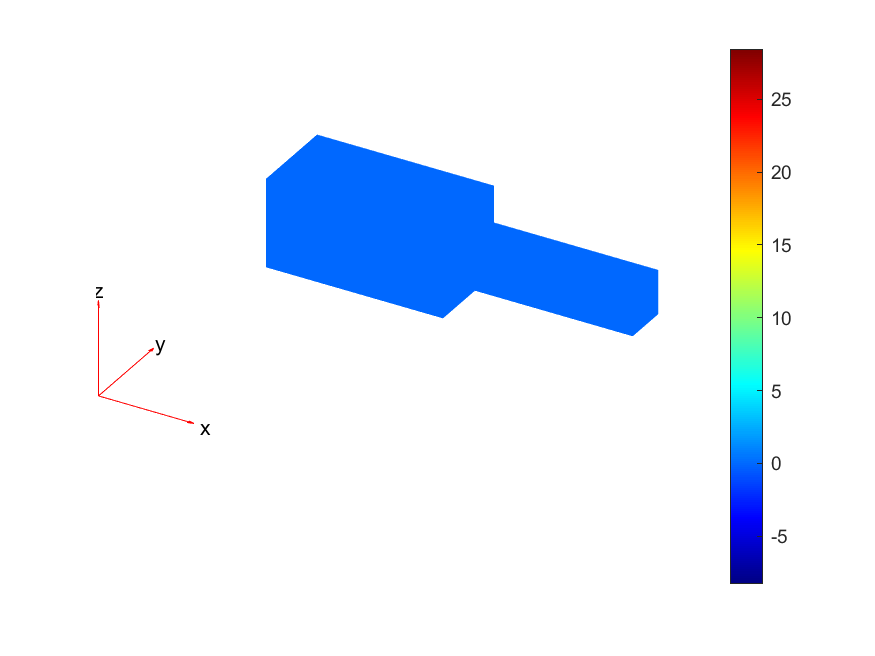}}
\scalebox{0.4}{\includegraphics{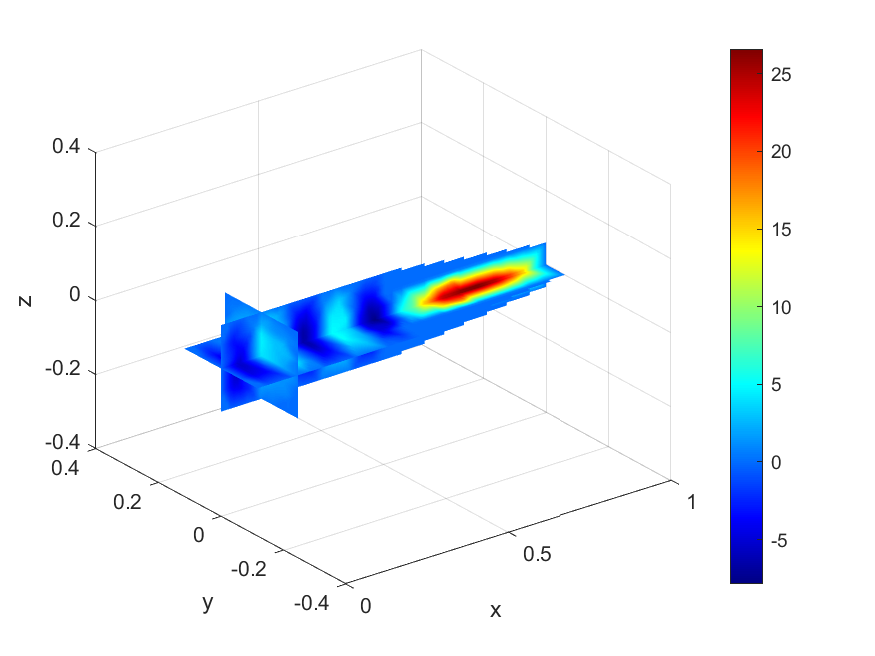}}
\end{center}
\caption{\small Approximations of eigenfunctions of \eqref{eq:ben26} with $G_\ee=(0,2^{-1})\times(-\ee,\ee)\times(-\ee,\ee) \, \cup \,
\{2^{-1}\}\times (-\ee,\ee)\times(-\ee,\ee) \, \cup \, (1/2,1)\times (-\ee2^{-1},\ee2^{-1})\times (-\ee2^{-1},\ee2^{-1})$ and $\ee=8^{-1}$. The figures are obtained choosing the eigenvalues $\lambda_1^\ee \approx 357.38,$ $\lambda_7^\ee\approx 1000.07$ and $\l_{12}^\ee\approx 1338.21.$} \label{fig:ben9}
\end{figure}

Figure \ref{fig:ben9}  shows numerical approximations of the eigenfunctions corresponding to the first,  seventh and twelfth eigenvalue of \eqref{eq:ben26} when the domain $G_\ee=(0,2^{-1})\times(-\ee,\ee)\times(-\ee,\ee) \, \cup \,
\{2^{-1}\}\times (-\ee,\ee)\times(-\ee,\ee) \, \cup \, (1/2,1)\times (-\ee2^{-1},\ee2^{-1}) $ $\times (-\ee2^{-1},\ee2^{-1})$ and $\ee=8^{-1}$ (see \eqref{eq:ben11}).
Here, $\lambda_1^\ee \approx 357.38,$ $\lambda_7^\ee\approx 1000.07$ and $\l_{12}^\ee\approx 1338.21.$
Also notice phenomena of localization.

\begin{figure}[t]
\begin{center}
\scalebox{0.4}{\includegraphics{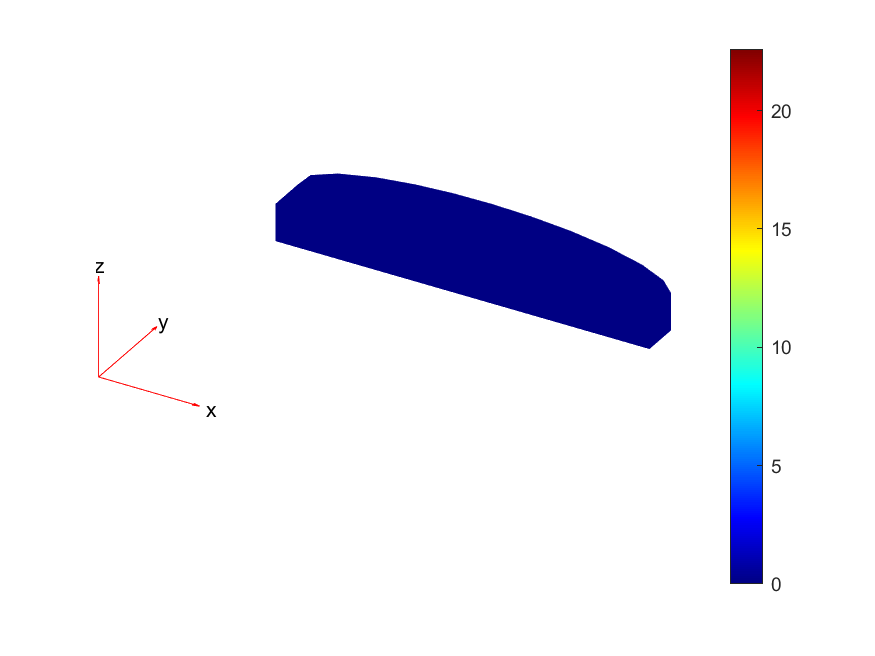}}
\scalebox{0.4}{\includegraphics{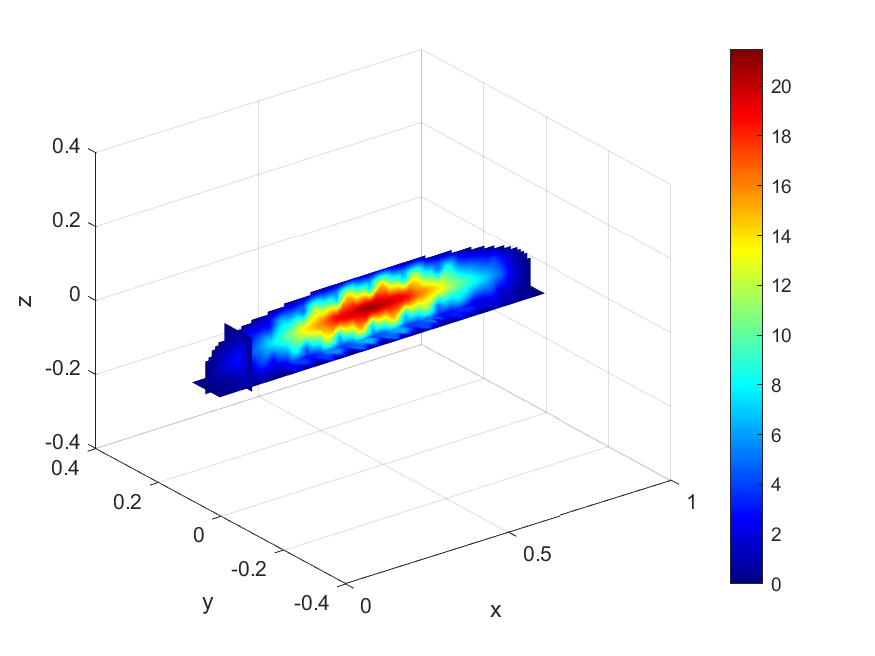}}
\scalebox{0.4}{\includegraphics{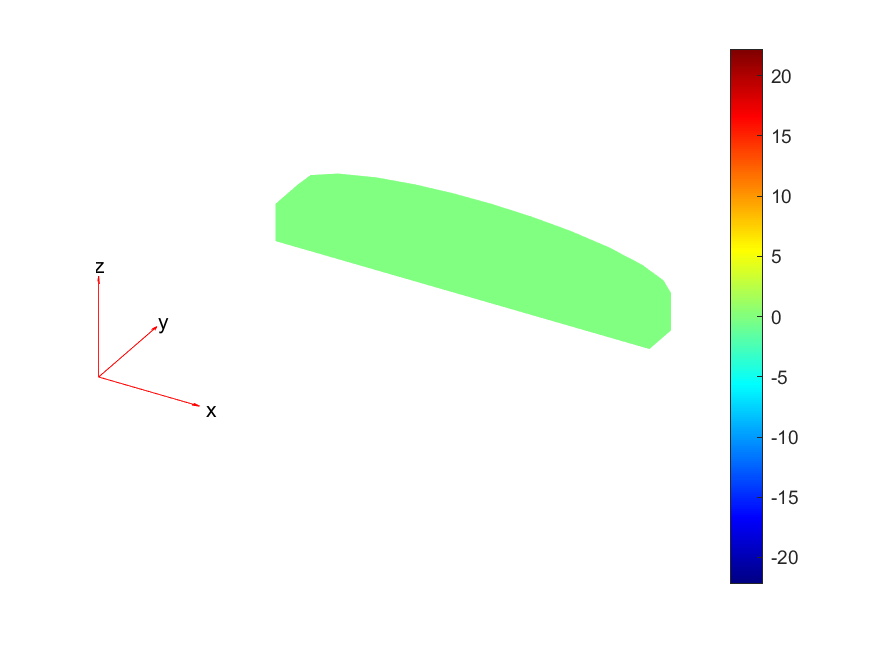}}
\scalebox{0.4}{\includegraphics{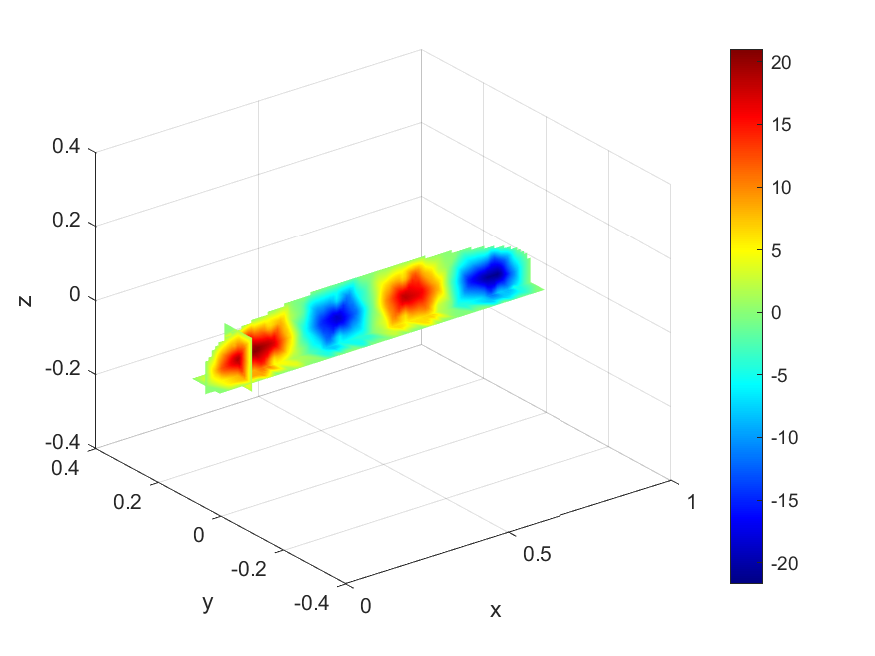}}
\scalebox{0.4}{\includegraphics{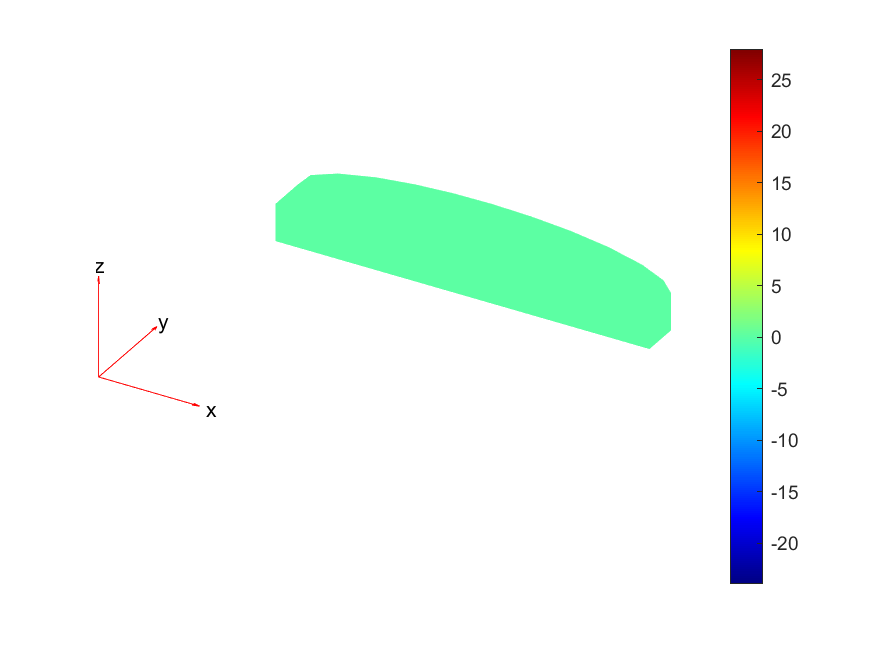}}
\scalebox{0.4}{\includegraphics{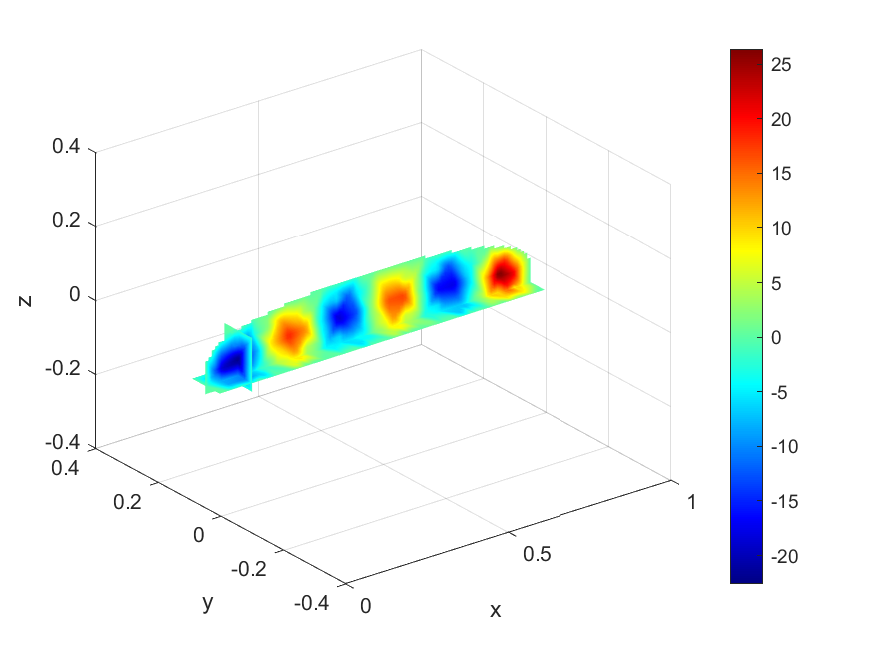}}
\end{center}
\caption{\small Approximations of eigenfunctions of \eqref{eq:ben26} with $G_\ee=(0,1)\times (0,\ee)\times (-\ee,\ee h(x_1))$ for a  certain function $h$ satisfying \eqref{eq:ben13} when $\ee=0.1$. The figures are obtained choosing the eigenvalues $\lambda_1^\ee\approx 1350.45,$ $\lambda_4^\ee\approx 1611.74$, $\lambda_6^\ee\approx 1925.74$. } \label{fig:ben10}
\end{figure}

Figure \ref{fig:ben10} shows numerical approximations of the eigenfunctions corresponding to the first, fourth and   sixth  eigenvalue of \eqref{eq:ben26} when the domain $G_\ee=(0,1)\times (0,\ee)\times (-\ee,\ee h(x_1))$ for a certain function $h$ satisfying \eqref{eq:ben13} and $\ee=0.1$ (see \eqref{eq:ben15}).
Here, $\lambda_1^\ee\approx 1350.45,$ $\lambda_4^\ee\approx 1611.74$ and  $\lambda_6^\ee\approx 1925.74$.

\begin{figure}[t]
\begin{center}
\scalebox{0.4}{\includegraphics{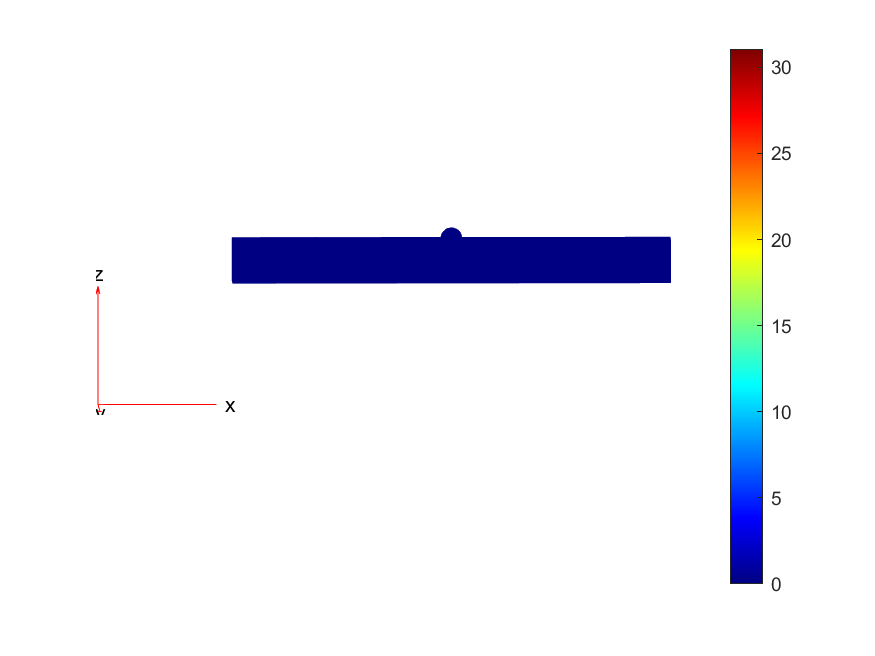}}
\scalebox{0.4}{\includegraphics{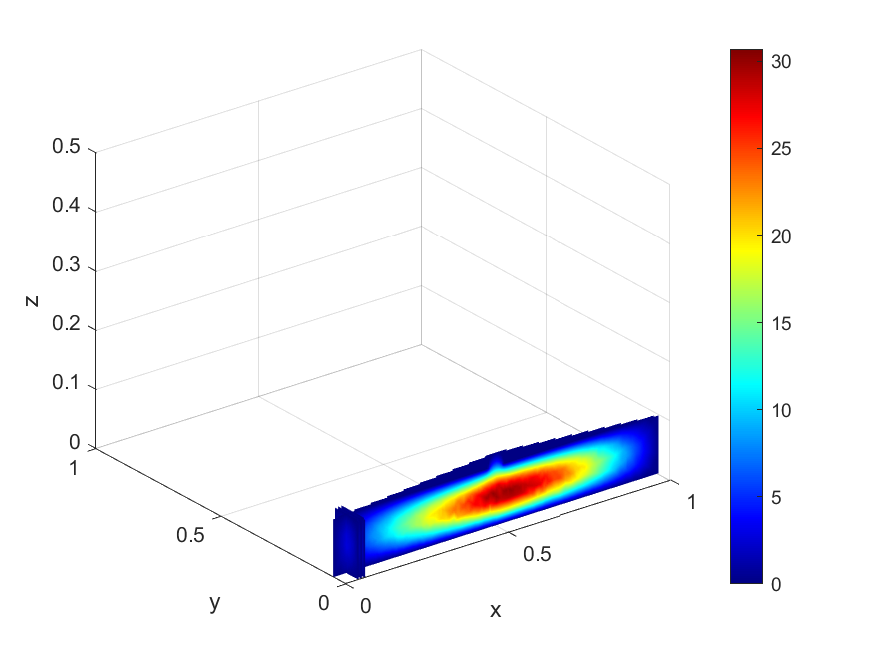}}
\end{center}
\caption{\small Approximation of the first eigenfunction of the Dirichlet problem  in a thin prism $(0,1)\times (0,\ee)\times (0,\ee)$ with a small perturbation when $\ee=0.1$. Here,  $\lambda_1^\ee\approx 2001.63.$  } \label{fig:ben1}
\end{figure}

Finally, notice that a small perturbation of the domain seems do not affect numerically to the first eigenfunction (see Figure \ref{fig:ben1} and compare with the two first graphics of Figure \ref{fig:ben8}).

\vspace{0.3cm}

\noindent
{\bf Acknowledgements} This work has been partially supported by grant PID2022-137694NB-I00 funded by MICIU/AEI/10.13039/501100011033 and by ERDF/EU.

\end{document}